%% file: document.tex
\documentclass[a4paper,11pt]{article}

\usepackage{xcolor}
\definecolor{navyblue}{RGB}{0,0,128}

\usepackage[left=3cm,right=3cm,top=2.5cm,bottom=3cm]{geometry}

\usepackage[colorlinks,
            linkcolor=navyblue,
            anchorcolor=black,
            citecolor=black,
            urlcolor=teal]{hyperref}

\usepackage{mathtools}
\usepackage{eucal}
\usepackage{amsmath,amssymb}
\usepackage{amsthm,thmtools}
\usepackage{bm}

\declaretheoremstyle[spaceabove=6pt,
                     spacebelow=9pt,
                     headfont=\normalfont\bfseries,
                     notefont=\mdseries,
                     notebraces={(}{)},
                     bodyfont=\normalfont,
                     numberwithin=section,
                     postheadspace=0.6em]{thmsty}
\declaretheoremstyle[spaceabove=6pt,
                     spacebelow=9pt,
                     headfont=\normalfont\bfseries,
                     notefont=\mdseries,
                     notebraces={(}{)},
                     bodyfont=\normalfont,
                     numberwithin=section,
                     postheadspace=0.6em,
                     sharenumber=theorem]{lemsty}
\declaretheoremstyle[spaceabove=6pt,
                     spacebelow=9pt,
                     headfont=\normalfont\bfseries,
                     bodyfont=\normalfont,
                     postheadspace=0.6em,
                     numbered=no,
                     qed=\qedsymbol]{prfsty}
\declaretheoremstyle[spaceabove=6pt,
                     spacebelow=9pt,
                     headfont=\normalfont\bfseries,
                     notefont=\mdseries,
                     notebraces={(}{)},
                     bodyfont=\normalfont,
                     numberwithin=section,
                     postheadspace=0.6em,
                     sharenumber=theorem]{defsty}
\declaretheoremstyle[spaceabove=6pt,
                     spacebelow=9pt,
                     headfont=\normalfont\bfseries,
                     bodyfont=\normalfont,
                     postheadspace=0.6em,
                     numberwithin=section,
                     sharenumber=theorem]{remsty}
\declaretheoremstyle[spaceabove=6pt,
                     spacebelow=9pt,
                     headfont=\normalfont\bfseries,
                     notefont=\mdseries,
                     notebraces={(}{)},
                     bodyfont=\normalfont,
                     postheadspace=0.6em,
                     numberwithin=section]{exasty}
\declaretheoremstyle[spaceabove=6pt,
                     spacebelow=9pt,
                     headfont=\normalfont\bfseries,
                     bodyfont=\normalfont,
                     numbered=no,
                     postheadspace=0.6em]{dittosty}

\declaretheorem[style=thmsty]{theorem}
\declaretheorem[style=thmsty]{lemma}
\declaretheorem[style=prfsty]{Proof}
\declaretheorem[style=thmsty]{definition}
\declaretheorem[style=thmsty]{remark}

\usepackage{enumerate}
\usepackage[toc,page]{appendix}
\usepackage{standalone,tikz}

\numberwithin{equation}{section} 
\DeclareMathOperator{\dif}{\mathrm{d}\!}
\DeclareMathOperator{\Det}{det}
\DeclareMathOperator{\Rank}{rank}
\renewcommand{\Re}{\operatorname{Re}}
\renewcommand{\Im}{\operatorname{Im}}
\newcommand{~}{\hskip0.065cm}

\title{Time-harmonic scattering of plane waves from \\ an infinite periodically inhomogeneous medium}
\author{Guanghui Hu, Andreas Rathsfeld, Jiayi Zhang, Ruming Zhang}
\date{\today}

\begin{document}
\maketitle	

\begin{abstract}
	We propose a new radiation condition for an infinite inhomogeneous two-dimensional 
    medium which is periodic in the vertical direction and remains invariant in 
    the horizontal direction. The classical Rayleigh-expansion radiation condition 
    does not apply to our case, because this would require the medium to be inhomogeneous 
	in a half plane. We utilize the Floquet theory to derive upward/downward wave 
	modes and define radiation conditions by expansions w.r.t.~these modes. The 
	downward radiation conditions leads to a downward Dirichlet-to-Neumann map 
	which can be used to truncate the infinite inhomogeneous domain in the vertical 
	direction. So we prove mapping properties of the upward/downward Dirichlet-to-Neumann 
	maps based on the asymptotic behavior of high-order wave modes. Finally, we 
    verify the strong ellipticity of the sesquilinear form corresponding to the 
	new scattering problem and show the unique solvability for all wavenumbers 
	with the exception of a countable set of numbers bounded below by a small 
	positive constant.
 \\

\end{abstract}


\section{Introduction}
Time-harmonic scattering problems for periodic gratings have attracted considerable
attention in the mathematical community (cf.~e.g.~\cite{Ba1,Ba2,Do,Ki}). A key issue 
for the scattering problems in unbounded domains is how to characterize an outgoing 
radiation condition at infinity. A physical-meaningful radiation condition should 
guarantee well-posedness of the mathematical model. Another motivation for studying 
radiation conditions is the design of efficient numerical schemes and the convergence 
analysis. With the help of radiation conditions, one can construct an appropriate 
Dirichlet-to-Neumann map to reduce the unbounded physical domain to a bounded 
computational domain. In periodic structures, the classical Rayleigh-expansion 
radiation condition applies to scattering of plane waves in homogeneous media. 
However, inhomogeneous media play an important role in practical applications, 
such as photonic crystals \cite{JoJoWiMe} and periodic waveguides, etc. The purpose 
of this work is to derive radiation conditions in an open periodic waveguide 
lying above an infinitely inhomogeneous medium and to prove solvability results 
for the scattering of plane waves by such kind of periodic media.

We consider the time-harmonic scattering of plane waves by a two-dimensional and 
horizontally periodic grating structure between a homogeneous cover and a substrate, 
which is periodically inhomogeneous w.r.t.~the vertical direction (see Figure \ref{fig:bvp}). 
Plane wave incidence together with the periodicity of the medium in horizontal 
direction gives rise to the quasi-periodicity of the wave fields. A challenging 
problem is to handle the infinity of the inhomogeneous medium in the substrate. 
Similarly to the generalized plane-wave modes in the derivation of the classical 
Rayleigh expansion, we can obtain wave modes by solving an ordinary differential 
equation with periodic coefficients. Then we can define the upward and downward 
radiation conditions by choosing outgoing bounded wave modes. Using these radiation 
conditions, we construct a downward Dirichlet-to-Neumann (DtN) map and formulate 
the boundary value problem in a periodic cell. We show some mapping properties 
of the DtN map and prove that the sesquilinear form in the variational formulation, 
corresponding to our boundary value problem, is strongly elliptic. Furthermore, 
uniqueness and existence of weak solutions are proved for all incident angles 
if the positive wavenumber is sufficiently small or if it is not in a countable 
set of exceptional wavenumbers with the only accumulating point at infinity. If 
the refractive index in the substrate depends only on the horizontal coordinate, 
upward and downward radiating modes can be derived by solving a one-dimensional 
Sturm-Liouville eigenvalue problem for an ODE in matrix form (see the authors' 
previous work \cite{HR2024}). This current work together with \cite{HR2024} provides 
a first insight into the well-posedness of scattering problems in an infinitely 
bi-periodic inhomogeneous medium. 
 
In closed periodic waveguides \cite{FlJo}, Fliss and Joly use the Floquet-Bloch 
theory and dispersion relations to classify the wave modes propagating into two 
different directions, under the assumption of non-vanishing group velocity for 
dispersion curves. Then outgoing radiation conditions are obtained by applying 
the Limiting Absorption Principle (LAP). Comparing \cite{FlJo} with the current 
work, we can define an explicit DtN operator without using LAP and with more mapping 
properties in Sobolve spaces. One reason for this lies in the fact that, due to 
plane-wave incidence, the solution modes depend explicitly on the the horizontal 
coordinate. In a series of works by Kirsch and coauthors \cite{K18, K22, K24}, 
radiation conditions for open periodic waveguides with local perturbations are 
derived mostly motivated by the LAP arguments of \cite{FlJo}. However, these are 
confined to homogeneous background media. Further, uniqueness and existence of 
solutions caused by a compactly supported source term are verified in those papers. 
In contrast, here we consider open periodic waveguides, but restrict our attention 
to scattering problems for plane-wave incidence into an infinite and inhomogeneous 
medium, without using LAP arguments. The mathematical setting of this work is 
close to that of Lamacz and Schweizer \cite{LaSc}, who derive a radiation condition 
for photonic crystals in a bi-periodic medium. Based on Bloch expansions and the 
Poynting vector, they construct an outgoing wave condition which only contains 
outgoing Bloch waves. Altogether, all of the above works extend the Rayleigh-expansion 
radiation condition from infinite homogeneous periodic settings to more complicated 
periodic materials. The advantage of our research methodology lies in its ability 
to obtain explicit wave modes and prove stronger solvability results by the DtN method.

This paper is organized as follows. In Section \ref{sec:problem}, we introduce
the medium, which is $2\pi$-periodic in vertical direction. We use the Floquet
theory to get the wave modes and then define the radiation condition in Section
\ref{sec:radiation}. In Section \ref{sec:solve}, we construct the Dirichlet-to-Neumann
map based on the new radiation condition. Then a boundary value problem for a
grating between a homogeneous and an inhomogeneous periodic medium is formulated.
We apply the variational method to discuss the solvability of the problem. Under
a few technical conditions, we show the existence of a unique solution for the
boundary value problem for any wavenumber $k$ not contained in an at most countable
set of exceptional wavenumbers. In other words, consider a horizontically periodic
grating structure above a substrate, the refractive index of which is periodic
in vertical and constant in horizontal direction. Then we prove that (cf.~Theorem \ref{mainth}),
with the exception of countable wavenumbers, the problem of plane waves scattered
by the grating has a unique solution.


\section{Quasiperiodic boundary value problem in an inhomogeneous half space} \label{sec:problem}

\begin{figure}[htp]
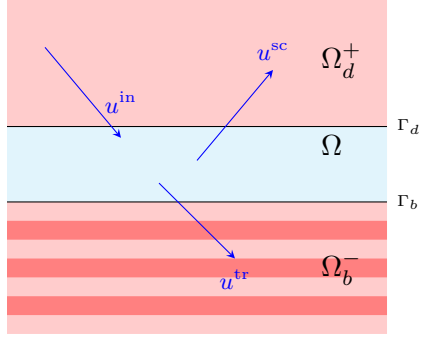

	\centering
	\includestandalone{figs/bvp}
	\caption{Geometry of the diffraction problem.}
	\label{fig:bvp}
\end{figure}

We first introduce the geometry of the problem, as shown in Figure \ref{fig:bvp}.
Choose real numbers $b$ and $d$ with $d > b$ and the two lines, $\Gamma_{d}
\coloneqq \{x = (x_{1}, x_{2})^{\mathrm{T}} \in \mathbb{R}^{2} \colon x_{2} = d\}$
and $\Gamma_{b} \coloneqq \{x = (x_{1}, x_{2})^{\mathrm{T}} \in \mathbb{R}^{2}
\colon x_{2} = b\}$, which divide $\mathbb{R}^{2}$ into the three regions
\begin{align*}
	&\Omega_{d}^{+} \coloneqq \{(x_{1}, x_{2})^{\mathrm{T}} \in \mathbb{R}^{2} \colon x_{2} > d\}, \\
	&\Omega\hskip0.25cm \coloneqq \{(x_{1}, x_{2})^{\mathrm{T}} \in \mathbb{R}^{2} \colon b < x_{2} < d\}, \\
	&\Omega_{b}^{-} \coloneqq \{(x_{1}, x_{2})^{\mathrm{T}} \in \mathbb{R}^{2} \colon x_{2} < b\}.
\end{align*}
The medium in $\mathbb{R}^{2}$ is characterized by the refraction index
$\tilde{q}(x_{1}, x_{2}) \in L^{\infty}(\mathbb{R}^{2})$ such that
\begin{equation*}
    \tilde{q}(x_{1}, x_{2}) =
    \begin{cases}
        1, &\text{in } \Omega_{d}^{+}, \\
        q_{0}(x_{1}, x_{2}), &\text{in } \Omega, \\
        q(x_{2}), &\text{in } \Omega_{b}^{-}.
    \end{cases}
\end{equation*}
Here the first function $x=(x_1, x_2)^{\mathrm{T}} \mapsto q_{0}(x_{1}, x_{2})$ is 
supposed to be {\tt p}-periodic w.r.t.~$x_{1}$ (${\tt p}>0$) and bounded w.r.t.~$x_2$ 
in $\Omega$. The second function $x_2 \mapsto q(x_2)$ is real-valued and bounded 
w.r.t.~$x_{2}\in (-\infty, b]$. Furthermore, we assume $\Im q_{0}(x_{1},x_{2}) \geqslant 0$,
$\Re q_0(x_1,x_2)\!>\!0$ for all $x\in \Omega$ and $q(x_2)\!>\!0$ for all $x_{2} < b$.
Let $u^{\mathrm{in}}(x_{1}, x_{2}) \coloneqq e^{i\hat{\alpha} x_{1} - i\beta x_{2}}$ be
the field incident from the region $\Omega_{d}^{+}$, where $\hat{\alpha} \coloneqq
k\sin\theta$, $\beta \coloneqq k\cos\theta$, $k$ is the positive wavenumber,
and $\theta \in (-\pi/2, \pi/2)$ the angle of incidence. Since the medium in 
$\Omega_{d}^{+}$ is homogeneous, $u^{\mathrm{in}}$ satisfies the Helmholtz equation
\begin{equation*}
	\Delta u^{\mathrm{in}} + k^2 u^{\mathrm{in}} = 0 \text{ in } \Omega_{d}^{+}.
\end{equation*}
The infinite region $\Omega_{b}^{-}$ is supposed to be filled by an infinitely
periodic medium characterized by the $2\pi$-periodic refractive index $q(x_{2})$.
Then the total field $u$ can be written as
\begin{equation*}
	u =
	\begin{cases}
		u^{\mathrm{in}} + u^{\mathrm{sc}} &\text{ in } \Omega_{d}^{+},\\
        u                                 &\text{ in } \Omega, \\
		u^{\mathrm{tr}}                   &\text{ in } \Omega_{b}^{-},
	\end{cases}
\end{equation*}
where $u^{\mathrm{sc}}$ is the scattered and $u^{\mathrm{tr}}$ the transmitted 
field. The time-harmonic acoustic wave propagation in $\mathbb{R}^{2}$ is governed 
by the Helmholtz equation
\begin{equation*}
	\Delta u(x_{1}, x_{2}) + k^{2} \tilde{q}(x_{1}, x_{2}) u(x_{1}, x_{2}) = 0\ , (x_{1}, x_{2})^{\mathrm{T}} \in \mathbb{R}^{2},
\end{equation*}
where $u$ denotes the pressure of an acoustic wave or a transverse field component
of a polarized electromagnetic wave.

\begin{remark}
	A radiation condition needs to be imposed at infinity to ensure the well-posedness
	of the diffraction problem. We recall that, in the case $q(x_{2}) \equiv 1$ in
    $\Omega_{b}^{-}$, the solution $u$ satisfies the downward radiating (outgoing)
    classical Rayleigh-expansion radiation condition if $u$ admits a Rayleigh expansion
	\begin{equation} \label{eq:clurc}
		u(x) = \sum_{n \in \mathbb{Z}} u_{n} e^{i(\alpha_{n}x_{1} - \beta_{n}(x_{2} - b))},\ u_{n} \in \mathbb{C},\ x_{2} < b,
	\end{equation}
	where $i = \sqrt{-1}$ and
	\begin{equation} \label{eq:alphanBetan}
		\alpha_{n} \coloneqq \hat{\alpha} + {\frac{2\pi}{{\tt p}}\,}n, \ \beta_{n} \coloneqq
		\begin{cases}
        \sqrt{k^{2} - \alpha_{n}^{2}} &\text{if } \alpha_{n}^{2} \,\le\; k^{2}, \\
        i\sqrt{\alpha_{n}^{2} - k^{2}} &\text{if }\alpha_{n}^{2} > k^{2}.
		\end{cases}
	\end{equation}
	For $\Omega_{b}^{-}$, we shall define generalized upward and downward propagating 
	radiation conditions later. These new radiation conditions will generalize 
	the above Rayleigh-expansion radiation condition from the infinite homogeneous 
	medium to an infinite medium, the refractive index $q$ of which is periodic 
	in the vertical direction but remains invariant in the horizontal direction.
\end{remark}

\begin{definition}
    The wave function $u$ is called $\hat{\alpha}$-quasi-periodic in $x_{1}$ with the
    parameter $\hat{\alpha} \coloneqq k \sin\theta$, $\theta \in (-\pi / 2, \pi / 2)$,
    if $x_{1} \mapsto\, u(x_{1}, x_{2})e^{-i \hat{\alpha} x_{1}}$ is {\tt p}-periodic
    for any fixed $x_{2}$.
\end{definition}

Define the quasi-periodic Sobolev spaces on $\Omega_{d}^{+}$, $\Omega_{b}^{-}$,
and $\mathbb{R}$ by
\begin{align*}
    H_{\hat{\alpha}}^{1}(\Omega_{d}^{+}) &\coloneqq \{u \in H_{loc}^{1}(\Omega_{d}^{+}) : u \text{ is } \hat{\alpha}\text{-quasi-periodic in } x_{1}\}, \\
	H_{\hat{\alpha}}^{1}(\Omega_{b}^{-}) &\coloneqq \{u \in H_{loc}^{1}(\Omega_{b}^{-}) : u \text{ is } \hat{\alpha}\text{-quasi-periodic in } x_{1}\}, \\
	H_{\hat{\alpha}}^{1/2}(\mathbb{R}) &\coloneqq \{f \in H_{loc}^{1/2}(\mathbb{R}) : e^{-i \hat{\alpha} x_{1}}f(x_{1}) \text{ is } {\tt p} \text{-periodic in } x_{1} \}, \\
	H_{\hat{\alpha}}^{-1/2}(\mathbb{R}) &\coloneqq \{f \in H_{loc}^{-1/2}(\mathbb{R}) : e^{-i \hat{\alpha} x_{1}}f(x_{1}) \text{ is } {\tt p} \text{-periodic in } x_{1} \}.
\end{align*}
Consider the Dirichlet boundary value problem in an inhomogeneous half space: For 
given Dirichlet data $f\in H_{\hat{\alpha}}^{1/2}(\mathbb{R})$, seek a solution
$u\in H^1_{\hat{\alpha}} (\Omega_b^{-})$ such that
\begin{equation} \label{eq:bvp}
	\begin{cases}
    	\Delta u(x_{1}, x_{2}) + k^{2}q(x_{2})u(x_{1}, x_{2}) = 0,\;
    	\hskip0.6225cm\quad (x_{1}, x_{2})^{\mathrm{T}}\in\, \Omega_{b}^{-},\\
    	\hskip3.75cm u(x_1,b) = f (x_1),\quad x_1\in\mathbb{R}.
	\end{cases}
\end{equation}
To ensure the well-posedness of (\ref{eq:bvp}), $u$ is required to satisfy an 
appropriate downward propagating radiation condition to be discussed within this 
paper.

We may expand the Dirichlet data $u \big\vert_{\Gamma_{b}} = f$ into the Fourier series
\begin{equation*}
    f(x_{1}) = \sum_{n \in \mathbb{Z}} f_{n} e^{i \alpha_{n} x_{1}},
\end{equation*}
where the coefficients $f_{n} \in \mathbb{C}$ satisfy $\sum_{n \in \mathbb{Z}}
(1+\vert n\vert^{2})^{1/2}\vert f_{n}\vert^{2} < \infty$. Indeed, since $u$ is 
$\alpha$-quasi-periodic, it admits the Fourier expansion
\begin{equation*}
	e^{-i \hat{\alpha} x_{1}} u(x_{1}, x_{2}) = \sum_{n \in \mathbb{Z}} u_{n}(x_{2})e^{i {\frac{2\pi}{\tt p}} nx_{1}},\ x_{2} < b,
\end{equation*}
or equivalently,
\begin{equation} \label{eq:expansion}
	u(x_{1}, x_{2}) = \sum_{n \in \mathbb{Z}} u_{n}(x_{2}) e^{i \alpha_{n} x_{1}},\ x_{2} < b.
\end{equation}
Inserting (\ref{eq:expansion}) into the Helmholtz equation, we find that
\begin{equation} \label{eq:1}
	\sum_{n \in \mathbb{Z}} \big[u_{n}''(x_{2}) + (k^{2}q(x_{2}) - \alpha_{n}^{2}) u_{n}(x_{2})\big] e^{i \alpha_{n} x_{1}} = 0.
\end{equation}
Below we shall confine our studies to a real-valued index function $q(x_{2})$ 
for $x_{2} < b$. By (\ref{eq:1}), we need to look for solutions $u_{n}$ to the 
Hill's equations
\begin{equation*}
	u_{n}''(x_{2}) + (k^{2}q(x_{2}) - \alpha_{n}^{2}) u_{n}(x_{2}) = 0,\ x_{2} < b,\ n \in \mathbb{Z},
\end{equation*}
for which the functions $(x_1,x_2)\mapsto u_n(x_2) e^{i \alpha_{n} x_{1}}$ with 
$n \in \mathbb{Z}$ physically represent downward and upward wave modes. We will 
use the Floquet theory (cf.~Appendix \ref{appendix:floquet}) to find an outgoing 
solution to the above Hill's equation in $x_{2} < b$ for every fixed $n \in \mathbb{Z}$, 
which represents the downward propagating wave mode as $x_{2} \rightarrow \infty$.


\section{Downward radiation condition for $q = q\left(x_{2}\right)$ in $x_2<b$.} \label{sec:radiation}

\subsection{Downward and upward radiation conditions in $\mathbb{R}^2$.}

In this section we extend $q(x_{2})$ from $x_{2} < b$ to $x_{2} \in \mathbb{R}$
by the $2\pi$-periodicity extension. Consider the one-dimensional ordinary
differential equations
\begin{equation} \label{eq:hill}
    u_{n}''(x_{2}) + (k^{2}q(x_{2}) - \alpha_{n}^{2}) u_{n}(x_{2}) = 0,\ x_{2} \in \mathbb{R},\ n \in \mathbb{Z}.
\end{equation}
We need to find two linearly independent solutions to \eqref{eq:hill}, in order
to classify the downward and upward wave modes for each $n\in \mathbb{Z}$. Let
$w_{n, 1}$ and $w_{n, 2}$ be the linearly independent solutions of (\ref{eq:hill}),
which satisfy the ``initial'' value conditions
\begin{equation} \label{eq:initial}
	\begin{aligned}
		w_{n, 1}(0) = 1,\ w_{n, 1}'(0) = 0; \\
		w_{n, 2}(0) = 0,\ w_{n, 2}'(0) = 1.
	\end{aligned}
\end{equation}
Consider the matrix
\begin{equation*}
	W_{n} \coloneqq
	\begin{pmatrix}
		w_{n, 1}(2\pi)  & w_{n, 2}(2\pi) \\
		w_{n, 1}'(2\pi) & w_{n, 2}'(2\pi)
	\end{pmatrix}	
	\in \mathbb{C}^{2 \times 2}.
\end{equation*}
Using Liouville's formula for the Wronskian \cite[Lemma 3.11]{Te}, we have
\begin{equation} \label{eq:detwrongski}
	\Det W_{n} = 1.
\end{equation}

\begin{definition} \label{def:characteristic}
	The eigenvalues $\lambda_{n, 1}$, $\lambda_{n, 2}$ of $W_{n}$ are called the
	characteristic multipliers of (\ref{eq:hill}). The numbers $\mu_{n, 1}$,
	$\mu_{n, 2}$ satisfying $-1/2 < \Im \mu_{n,j}  \leqslant 1/2$ and $e^{2\pi\mu_{n, j}} = \lambda_{n, j}$,
	$j = 1, 2$, are called the characteristic exponents of (\ref{eq:hill}).
\end{definition}

It follows from (\ref{eq:detwrongski}) that
\begin{align}
	&\lambda_{n, j}^{2} - \eta_{n} \lambda_{n, j} + 1 = 0,\ j = 1, 2,\ n \in \mathbb{Z}, \label{eq:polynomial} \\
	&\eta_{n} \coloneqq w_{n, 1}(2\pi) + w_{n, 2}'(2\pi) \in \mathbb{C}. \label{eq:etan}
\end{align}
Consequently, we have
\begin{lemma}
	The characteristic multipliers and exponents of (\ref{eq:hill}) satisfy the
	relations
	\begin{equation*}
		\lambda_{n, 1} \lambda_{n, 2} = 1,\ \mu_{n, 1} + \mu_{n, 2} \in \{0, i\},\ \text{for all } n \in \mathbb{Z}.
	\end{equation*}
\end{lemma}

Without loss of generality we suppose that $\vert \lambda_{n, 1}\vert \geqslant
\vert \lambda_{n, 2}\vert$. Applying the Floquet theory we shall derive two linearly
independent solutions of (\ref{eq:hill}).

Now suppose $q(x_{2})$ is real-valued. Then the $w_{n, j}$ and $\eta_{n}$ are 
all real-valued.

\begin{lemma} \label{lem:solvereal}
  There are two linearly independent solutions $\xi_{n, j}$, $j = 1, 2$ of (\ref{eq:hill})
  expressed with the help of two $2\pi$-periodic functions $v_{n, j}$, $j = 1, 2$,
  where these expressions depend on the number $\eta_{n}$ as follows.
	\begin{enumerate}[(a)]
		\item If $\eta_{n} > 2$, then
			  \begin{equation*}
			  	\xi_{n, 1}(x_{2}) = e^{\mu_{n, 1}x_{2}}v_{n, 1}(x_{2}),\ \xi_{n, 2}(x_{2}) = e^{-\mu_{n, 1}x_{2}}v_{n, 2}(x_{2}),
			  \end{equation*}
		  	  where $\mu_{n, 1}$ is a positive real number such that $e^{2\pi\mu_{n, 1}} = \lambda_{n, 1}$.
		\item If $\eta_{n} < -2$, then
		      \begin{equation*}
		          \xi_{n, 1}(x_{2}) = e^{\left(\tilde{\mu}_{n, 1} + \frac{i}{2}\right) x_{2}}v_{n, 1}(x_{2}),\ \xi_{n, 2}(x_{2}) = e^{\left(-\tilde{\mu}_{n, 1} + \frac{i}{2}\right) x_{2}} v_{n, 2}(x_{2}),
		      \end{equation*}
	      	  where $\tilde{\mu}_{n, 1}$ is a positive real number such that $e^{2\pi \tilde{\mu}_{n, 1}} = \vert \lambda_{n, 1}\vert$.
		\item If $\vert \eta_{n}\vert < 2$, then
		      \begin{equation*}
		      	\xi_{n, 1}(x_{2}) = e^{-i \theta_{n} x_{2}}v_{n, 1}(x_{2}),\ \xi_{n, 2}(x_{2}) = e^{i \theta_{n} x_{2}}v_{n, 2}(x_{2}),
		      \end{equation*}
	          where $\theta_{n}=\mu_{n,1}$ is a real number and $\theta_{n} \in \left(0, 1/2
	          \right)$ such that $e^{i 2\pi\theta_{n}} = \lambda_{n, 1}$.
		\item If $\eta_{n} = 2$, then there are two sub-cases as follows:
		      \begin{enumerate}[\text{case} (i):]
		      	\item For $\Rank(W_{n} - I) = 0$, we have
                      \begin{equation*}
		      		  	   \xi_{n, 1}(x_{2}) = v_{n, 1}(x_{2}),\ \xi_{n, 2}(x_{2}) = v_{n, 2}(x_{2}),
		      	      \end{equation*}
	      		\item For $\Rank(W_{n} - I) \neq 0$, we have
                      \begin{equation*}
	      			       \xi_{n, 1}(x_{2}) = v_{n, 1}(x_{2}),\ \xi_{n, 2}(x_{2}) = x_{2}v_{n, 1}(x_{2}) + v_{n, 2}(x_{2}).
	      		      \end{equation*}
		      \end{enumerate}
		\item If $\eta_{n} = -2$, then there are two sub-cases:
		      \begin{enumerate}[\text{case} (i):]
		      	\item For $\Rank(W_{n} + I) = 0$, we have
                      \begin{equation*}
		      		       \xi_{n, 1}(x_{2}) = e^{i \frac{x_{2}}{2}} v_{n, 1}(x_{2}),\ \xi_{n, 2}(x_{2}) = e^{i \frac{x_{2}}{2}} v_{n, 2}(x_{2}),
		      	      \end{equation*}
		      	\item For $\Rank(W_{n} + I) \neq 0$, we have
                      \begin{equation*}
		      		       \xi_{n, 1}(x_{2}) = e^{i \frac{x_{2}}{2}} v_{n, 1}(x_{2}),\ \xi_{n, 2}(x_{2}) = e^{i \frac{x_{2}}{2}} (x_{2}v_{n, 1}(x_{2}) + v_{n, 2}(x_{2})).
		      	      \end{equation*}
		      \end{enumerate}
	\end{enumerate}
\end{lemma}

\begin{Proof}
	\begin{enumerate}[(a)]
		\item Since $\eta_{n} > 2$ and $\lambda_{n, 1}\lambda_{n,2} = 1$, $\lambda_{n,
              1}$ and $\lambda_{n, 2}$ are distinct positive real numbers and 
			  $0 < \lambda_{n, 2} < 1 < \lambda_{n, 1}$. By Definition \ref{def:characteristic},
              there exists a positive real number $\mu_{n, 1}$ such that
		      \begin{equation*}
			       e^{2\pi \mu_{n, 1}} = \lambda_{n, 1},\ e^{-2\pi \mu_{n, 1}} = \lambda_{n, 2}.
		      \end{equation*}
		      Thus, by part (1) of Theorem \ref{thm:floquet},
		      \begin{equation*}
			       \xi_{n, 1}(x_{2}) = e^{\mu_{n, 1}x_{2}}v_{n, 1}(x_{2}),\ \xi_{n, 2}(x_{2}) = e^{-\mu_{n, 1}x_{2}}v_{n, 2}(x_{2}),
		      \end{equation*}
	          where $v_{n, 1}$ and $v_{n, 2}$ are $2\pi$-periodic functions.
	    \item Since $\eta_{n} < -2$ and $\lambda_{n, 1}\lambda_{n, 2} = 1$, $\lambda_{n,
              1}$ and $\lambda_{n, 2}$ are both negative real numbers and $\lambda_{n, 1}
              < -1 < \lambda_{n, 2} < 0$. There exists a positive real number $\tilde{\mu}_{n, 1}$ such that
	          \begin{equation*}
	    	       e^{2\pi \tilde{\mu}_{n, 1}} = \vert \lambda_{n, 1}\vert,\ e^{-2\pi \tilde{\mu}_{n, 1}} = \vert \lambda_{n, 2}\vert.
	          \end{equation*}
       	      By part (1) of Theorem \ref{thm:floquet}, we get
       		  \begin{equation*}
       			   \xi_{n, 1}(x_{2}) = e^{\left(\tilde{\mu}_{n, 1} + \frac{i}{2}\right) x_{2}}v_{n, 1}(x_{2}),\ \xi_{n, 2}(x_{2}) = e^{\left(-\tilde{\mu}_{n, 1} + \frac{i}{2}\right) x_{2}}v_{n, 2}(x_{2}),
       		  \end{equation*}
              where $v_{n, 1}$ and $v_{n, 2}$ are $2\pi$-periodic functions.
        \item If $\vert \eta_{n}\vert < 2$, then $\lambda_{n, 1}$ and $\lambda_{n, 2}$
		      are complex conjugates. Additionally, $\lambda_{n, 1}\lambda_{n, 2} = \vert \lambda_{n, 1}
		      \vert^{2} = 1$. We assume that $\lambda_{n, 1}$ lies in the upper half plane.
		      There is a real number $\theta_{n} \in (0, 1/2)$ such that
              \begin{equation*}
        	       e^{i2\pi \theta_{n}} = \lambda_{n, 1},\ e^{-i2\pi \theta_{n}} = \lambda_{n, 2}.
              \end{equation*}
    	      Therefore, by part (1) of Theorem \ref{thm:floquet},
    	      \begin{equation*}
    		      \xi_{n, 1}(x_{2}) = e^{-i \theta_{n} x_{2}}v_{n, 1}(x_{2}),\ \xi_{n, 2}(x_{2}) = e^{i \theta_{n} x_{2}}v_{n, 2}(x_{2}), 
    	      \end{equation*}
              where $v_{n, 1}$ and $v_{n, 2}$ are $2\pi$-periodic functions.
        \item If $\eta_{n} = 2$, then we have $\lambda_{n, 1} = \lambda_{n, 2} = 1$.
		      Under this condition, we need to compute the rank of $W_{n} - I \in \mathbb{C}^{2\times
              2}$ to decide which part of Theorem \ref{thm:floquet} is applicable.
              \begin{enumerate}[\text{case} (i):]
        	      \item $w_{n, 2}(2\pi) = w_{n, 1}'(2\pi) = 0$.

        	            Then
        	            \begin{equation*}
        	      	       \Det \begin{pmatrix}
        	      		    w_{n, 1}(2\pi) & w_{n, 1}'(2\pi) \\
        	      		    w_{n, 2}(2\pi) & w_{n, 2}'(2\pi)
        	      	      \end{pmatrix} =
              	          \Det \begin{pmatrix}
              	    	    w_{n, 1}(0) & w_{n, 1}'(0) \\
              	    	    w_{n, 2}(0) & w_{n, 2}'(0)
              	          \end{pmatrix} = 1.         	
        	            \end{equation*}
                        We have $w_{n, 1}(2\pi)w_{n, 2}'(2\pi) = 1$ and
                        \begin{equation*}
                  	       \eta_{n} = w_{n, 1}(2\pi) + w_{n, 2}'(2\pi) = 2.
                        \end{equation*}
              	        Hence $w_{n, 1}(2\pi) = w_{n, 2}'(2\pi) = 1$, which implies
              	        $\Rank\left(W_{n} - I\right) = 0$, i.e., there are two linearly
              	        independent solutions of (\ref{eq:hill}). Now we may apply
              	        part (1) of Theorem \ref{thm:floquet}. Since $\lambda_{n, 1}
              	        = \lambda_{n, 2} = 1$, we get $\mu_{n, 1} = \mu_{n, 2} = 0$ and
              	        \begin{equation*}
              	  	      \xi_{n, 1}(x_{2}) = v_{n, 1}(x_{2}),\ \xi_{n, 2}(x_{2}) = v_{n, 2}(x_{2}),
              	        \end{equation*}
                        where $v_{n, 1}$ and $v_{n, 2}$ are $2\pi$-periodic functions.
        	      \item $w_{n, 2}(2\pi), w_{n, 1}'(2\pi)$ not both zero
        	
        		        In this case $\Rank (W_{n} - I) \neq 0$. Applying part
        		        (2) of Theorem \ref{thm:floquet} gives $\mu_{n, 1} = \mu_{n, 2} = 0$ and
        		        \begin{equation*}
        		   	       \xi_{n, 1}(x_{2}) = v_{n, 1}(x_{2}),\ \xi_{n, 2}(x_{2}) = x_{2}v_{n, 1}(x_{2}) + v_{n, 2}(x_{2}),
        		        \end{equation*}
        	            where $v_{n, 1}$ and $v_{n, 2}$ are $2\pi$-periodic functions.
              \end{enumerate}
        \item Similar to the assertion (d), we consider the following two cases.
              \begin{enumerate}[\text{case} (i):]
              	  \item $w_{n, 2}(2\pi) = w_{n, 1}'(2\pi) = 0$.
              	
              	        We have that $\Rank(W_{n} + I) = 0$. By part (1) of 
						Theorem \ref{thm:floquet}, one obtains $\mu_{n, 1} = \mu_{n, 2} = i/2$ and
              	        \begin{equation*}
              	            \xi_{n, 1}(x_{2}) = e^{i \frac{x_{2}}{2}} v_{n, 1}(x_{2}),\ \xi_{n, 2}(x_{2}) = e^{i \frac{x_{2}}{2}} v_{n, 2}(x_{2}),
              	        \end{equation*}
                        where $v_{n, 1}$ and $v_{n, 2}$ are $2\pi$-periodic functions.
              	  \item $w_{n, 2}(2\pi), w_{n, 1}'(2\pi)$ not both zero.
              	
              	  		Here $\Rank(W_{n} + I) \neq 0$. Again using part (2) of Theorem
                        \ref{thm:floquet}, we get $\mu_{n, 1} = \mu_{n, 2} = i/2$ and
              	  		\begin{equation*}
              	  			\xi_{n, 1}(x_{2}) = e^{i\frac{x_{2}}{2}}v_{n, 1}(x_{2}),\ \xi_{n, 2}(x_{2}) = e^{i \frac{x_{2}}{2}}(x_{2} v_{n, 1}(x_{2}) + v_{n, 2}(x_{2})),
              	  		\end{equation*}
                		where $v_{n, 1}$ and $v_{n, 2}$ are $2\pi$-periodic functions.
              \end{enumerate}
	\end{enumerate}
\end{Proof}
\begin{remark}
  In Lemma \ref{lem:solvereal}, the functions $\xi_{n, j}$ are periodic in
  cases (d) (i) and anti-periodic in cases (e) (i).
\end{remark}

According to Lemma \ref{lem:solvereal}, we can now define the upward modes
$u_{n}^{+}(x_{2})e^{i\alpha_{n}x_{1}}$ {(that is, modes decaying for $x_{2} \rightarrow +\infty$
or modes propagating upwards)} and downward modes $u_{n}^{-}(x_{2}) e^{i\alpha_{n}x_{1}}$
(modes decaying for $x_2\rightarrow -\infty$ or modes propagating downwards),
where $u_{n}^{\pm}(x_{2})$ are defined as follows.
\begin{definition} \label{def:modes}
	\begin{enumerate}[(a)]
		\item If $\eta_{n} > 2$, then
		      \begin{equation*}
			      u_{n}^{+}(x_{2}) \coloneqq e^{-\mu_{n, 1}x_{2}}v_{n, 2}(x_{2}),\ u_{n}^{-}(x_{2}) \coloneqq e^{\mu_{n, 1}x_{2}}v_{n, 1}(x_{2}),
		      \end{equation*}
              where $v_{n, 1}$ and $v_{n, 2}$ are the $2\pi$-periodic functions 
			  of Lemma \ref{lem:solvereal} (a).
		\item If $\eta_{n} < -2$, then
			  \begin{equation*}
				  u_{n}^{+}(x_{2}) \coloneqq e^{\left(-\tilde{\mu}_{n, 1} + \frac{i}{2}\right) x_{2}}v_{n, 2}(x_{2}),\
			      u_{n}^{-}(x_{2}) \coloneqq e^{\left(\tilde{\mu}_{n, 1} + \frac{i}{2}\right) x_{2}}v_{n, 1}(x_{2}),
		      \end{equation*}
		      where $v_{n, 1}$ and $v_{n, 2}$ are the $2\pi$-periodic functions 
			  of Lemma \ref{lem:solvereal} (b).
		\item If $\vert \eta_{n}\vert < 2$, then
		      \begin{equation*}
			      u_{n}^{+}(x_{2}) \coloneqq e^{i \theta_{n} x_{2}} v_{n, 2}(x_{2}),\
			      u_{n}^{-}(x_{2}) \coloneqq e^{-i \theta_{n} x_{2}} v_{n, 1}(x_{2}),
		      \end{equation*}
		      where $v_{n, 1}$ and $v_{n, 2}$ are the $2\pi$-periodic functions 
			  of Lemma \ref{lem:solvereal} (c).
		\item If $\eta_{n} = 2$, then
		      \begin{enumerate}[\text{case} (i):]
			      \item \begin{equation*}
                            u_{n}^{+}(x_{2}) = u_{n}^{-}(x_{2}) \coloneqq (v_{n, 1}(x_{2}), v_{n, 2}(x_{2}))^{\mathrm{T}},
                        \end{equation*}
			      \item \begin{equation*}
	                        u_{n}^{+}(x_{2}) = u_{n}^{-}(x_{2}) \coloneqq v_{n, 1}(x_{2}),
			            \end{equation*}
		      \end{enumerate}
		      where $v_{n, 1}$ and $v_{n, 2}$ are the $2\pi$-periodic functions 
			  of Lemma \ref{lem:solvereal} {(d)}.
		\item If $\eta_{n} = -2$, then
		      \begin{enumerate}[\text{case} (i):]
			      \item \begin{equation*}
				            u_{n}^{+}(x_{2}) = u_{n}^{-}(x_{2}) \coloneqq  e^{i\frac{x_{2}}{2}}\left(v_{n, 1}(x_{2}), v_{n, 2}(x_{2})\right)^{\mathrm{T}},
			            \end{equation*}
			      \item \begin{equation*}
				        	u_{n}^{+}(x_{2}) = u_{n}^{-}(x_{2}) \coloneqq e^{i\frac{x_{2}}{2}}v_{n, 1}(x_{2}),
			            \end{equation*}
		      \end{enumerate}
		      where $v_{n, 1}$ and $v_{n, 2}$ are the $2\pi$-periodic functions 
			  of Lemma \ref{lem:solvereal} {(e)}.
	\end{enumerate}
    Note that $u_{n}^{\pm}$ are vectorial functions in cases (d) (i) and (e) (i).
\end{definition}

With the help of Definition \ref{def:modes}, we can define the upward and downward
radiation conditions by a possible representation as a superposition of upward and
downward modes, respectively.
\begin{definition}[Radiation Conditions] \label{def:radiationCondition}
    Consider $\hat{\alpha}$-quasiperiodic solutions to the Helmholtz equation 
	$\Delta u(x_{1}, x_{2}) + k^{2}q(x_{2})u(x_{1}, x_{2}) = 0$. The solution 
	$u$ is called an upward radiating solution if there exist an $a^{+} \in \mathbb{R}$ 
	and coefficients $C_{n}^{+}$ such that
    \begin{align}
        u(x_{1}, x_{2}) = \sum_{n \in \mathbb{Z}} C_{n}^{+}\cdot u_{n}^{+}(x_{2})\; e^{i\alpha_{n}x_{1}}\quad \text{ in } x_{2} > a^{+}, \tag{URC} \label{eq:urc}
    \end{align}
    where $C_{n}^{+} \in \mathbb{C}^2$ in case (d) (i) or (e) (i) and $C_{n}^{+} 
	\in \mathbb{C}$ otherwise.

    The solution $u$ is called an downward radiating solution if there exist an  $a^{-} \in \mathbb{R}$ and coefficients  $C_{n}^{-}$ such that
    \begin{align}
    	u(x_{1}, x_{2}) = \sum_{n \in \mathbb{Z}} C_{n}^{-}\cdot u_{n}^{-}(x_{2})\; e^{i\alpha_{n}x_{1}} \quad \text{ in } x_{2} < a^{-}, \tag{DRC} \label{eq:drc}
    \end{align}
    where $C_{n}^{-} \in \mathbb{C}^2$ in case (d) (i) or (e) (i) and $C_{n}^{-} 
	\in \mathbb{C}$ otherwise.
\end{definition}

\subsection{Special instance: $q(x_{2}) \equiv 1$}

It is natural to ask how do the new radiation conditions of Definition \ref{def:radiationCondition} 
generalize  the classical Rayleigh-expansion radiation conditions. Below we show 
that the Rayleigh-expansion condition in the homogeneous case that $q(x_{2}) \equiv 1$
is a special instance of Definition \ref{def:radiationCondition}. Consider 
$u_{n}''(x_{2}) + (k^{2} - \alpha_{n}^{2})u_{n}(x_{2}) = 0$. Here, $\alpha_{n}$ 
and $\beta_{n}$ have been defined in (\ref{eq:alphanBetan}).
\begin{enumerate}[(1)]
    \item If $k^{2} > \alpha_{n}^{2}$, then there are two linearly independent
          solutions
		  \begin{equation*}
		       \xi_{1}(x_{2}) = e^{i\sqrt{k^{2} - \alpha_{n}^{2}}x_{2}},\ \xi_{2}(x_{2}) = e^{-i\sqrt{k^{2} - \alpha_{n}^{2}}x_{2}}.
		  \end{equation*}
	      This corresponds to part (1) of Theorem \ref{thm:floquet} with $m_{1} =
          i\sqrt{k^{2}-\alpha_{n}^{2}}$, $m_{2} = -i\sqrt{k^{2}-\alpha_{n}^{2}}$
	      and $v_{1} = v_{2} = 1$, that is, case (c) in Lemma \ref{lem:solvereal}.
	\item If $k^{2} < \alpha_{n}^{2}$, then there are two linearly independent
          solutions
          \begin{equation*}
              \xi_{1}(x_{2}) = e^{\sqrt{\alpha_{n}^{2}-k^{2}}x_{2}},\ \xi_{2}(x_{2}) = e^{-\sqrt{\alpha_{n}^{2}-k^{2}}x_{2}}.
          \end{equation*}
          This corresponds to part (1) of Theorem \ref{thm:floquet} with
          $m_{1} = \sqrt{\alpha_{n}^{2} - k^{2}}$, $m_{2} = -\sqrt{\alpha_{n}^{2}
          - k^{2}}$ and $v_{1} = v_{2} = 1$, that is, case (a) in Lemma \ref{lem:solvereal}.
	\item If $k^{2} = \alpha_{n}^{2}$, then there are two linearly independent solutions
		  \begin{equation*}
		      \xi_{1}(x_{2}) = a,\ \xi_{2}(x_{2}) = ax_{2}+b, \text{ for some } a, b \in \mathbb{C}.
		  \end{equation*}
	      This corresponds to part (2) of Theorem \ref{thm:floquet} with $m = 0$,
          $v_{1} = a$ and $v_{2} = b$, that is, case (ii) in Lemma \ref{lem:solvereal} (d).
\end{enumerate}

Furthermore, we can also compute $\eta_{n}$ for each case above. Solving two initial 
value problems for the ODE
\begin{equation} \label{eq:ode}
	\begin{cases}
		w_{n, 1}''(x_{2}) + \beta_{n}^{2} w_{n, 1}(x_{2}) = 0, \\
		w_{n, 1}(0) = 1,\ w_{n, 1}'(0) = 0,
	\end{cases}
	\quad
	\begin{cases}
		w_{n, 2}''(x_{2}) + \beta_{n}^{2} w_{n, 2}(x_{2}) = 0, \\
		w_{n, 2}(0) = 0,\ w_{n, 2}'(0) = 1,
	\end{cases}
\end{equation}
we get solutions for (\ref{eq:ode}) as below:
\begin{align*}
	w_{n, 1}(x_{2})  &= \frac{1}{2}\left(e^{i \beta_{n} x_{2}} + e^{-i \beta_{n} x_{2}}\right), \\
	w_{n, 1}'(x_{2}) &= \frac{i \beta_{n}}{2}\left(e^{i \beta_{n} x_{2}} - e^{-i \beta_{n} x_{2}}\right), \\
	w_{n, 2}(x_{2})  &= \frac{1}{2i \beta_{n}}\left(e^{i \beta_{n} x_{2}} - e^{- i \beta_{n} x_{2}}\right) = \frac{1}{(i \beta_{n})^{2}} w_{n, 1}'(x_{2}), \\
	w_{n, 2}'(x_{2}) &= \frac{1}{2} \left(e^{i \beta_{n} x_{2}} + e^{-i \beta_{n} x_{2}}\right) = w_{n, 1}(x_{2}).
\end{align*}
Then we can compute $\eta_{n}$ as
\begin{equation*}
	\eta_{n} = w_{n, 1}(2\pi) + w_{n, 2}'(2\pi) = e^{i 2\pi \beta_{n}} + e^{-i 2\pi \beta_{n}}.
\end{equation*}

\begin{enumerate}[(1)]
    \item If $k^{2} > \alpha_{n}^{2}$, then $\eta_{n} = 2\cos(2 \pi \beta_{n})$
          and $\vert \eta_{n}\vert < 2$. This corresponds to case (c) in Lemma
          \ref{lem:solvereal}. We can define the upward mode $u_{n}^{+}(x_{2})e^{i\alpha_{n}x_{1}}$
          and the downward mode $u_{n}^{-}(x_{2})e^{i\alpha_{n}x_{1}}$, where
          \begin{equation*}
		       u_{n}^{+}(x_{2}) = e^{i \beta_{n} x_{2}} \text{ and } u_{n}^{-}(x_{2}) = e^{-i \beta_{n} x_{2}}.
	      \end{equation*}
    \item If $k^{2} < \alpha_{n}^{2}$, then $\beta_{n} = i\vert \beta_{n}\vert$
          and $\eta_{n} = e^{-2\pi \vert \beta_{n}\vert} + e^{2\pi \vert
          \beta_{n}\vert} > 2$. This corresponds to case (a) in Lemma \ref{lem:solvereal}.
          We can define the upward mode $u_{n}^{+}(x_{2})e^{i\alpha_{n}x_{1}}$ and the downward mode $u_{n}^{-}(x_{2})e^{i\alpha_{n}x_{1}}$, where
          \begin{equation*}
		       u_{n}^{+}(x_{2}) = e^{-\vert \beta_{n}\vert x_{2}} \text{ and } u_{n}^{-}(x_{2}) = e^{\vert \beta_{n}\vert x_{2}}.
	      \end{equation*}
    \item If $k^{2} = \alpha_{n}^{2}$, then $\beta_{n} = 0$ and $\eta_{n} = 2$.
          Check the two initial value problems
	      \begin{equation*}
		       \begin{cases}
			       w_{n, 1}'' = 0, \\
			       w_{n, 1}(0) = 1,\ w_{n, 1}'(0) = 0,
		       \end{cases}
		       \quad
		       \begin{cases}
			       w_{n, 2}'' = 0, \\
			       w_{n, 2}(0) = 0,\ w_{n, 2}'(0) = 1.
		       \end{cases}
	       \end{equation*}
	      We have $w_{n, 1} = 1,\ w_{n, 2} = x_{2}$, and
	      \begin{equation*}
		    W_{n} = \begin{pmatrix}
			        w_{n, 1}(2\pi) & w_{n, 1}'(2\pi) \\
			        w_{n, 2}(2\pi) & w_{n, 2}'(2\pi)
		    \end{pmatrix} =
		    \begin{pmatrix}
			    1    & 0 \\
			    2\pi & 1
		    \end{pmatrix},
	      \end{equation*}
	      $\Rank(W_{n} - I) = 1$. This corresponds to part (ii) of case (d) in 
		  Lemma \ref{lem:solvereal} with $v_{n,1}(x_2)\equiv 1$. We can only 
		  define the outgoing bounded wave mode $u_{n}^{\pm}(x) \coloneqq e^{i\alpha_{n}x_{1}}$.
\end{enumerate}

\subsection{Asymptotics of $\eta_{n}$ and $\mu_{n}$ as $\vert n\vert \to +\infty$.}

The purpose of this subsection is to prove the asymptotic behavior of $\eta_n$ 
and $\mu_n$ as $|n| \rightarrow \infty$ through the comparison theorem for initial 
value problems of an ordinary differential equation, which is stated as below.

\begin{lemma}[Comparison Theorem] \label{lem:compareOde}
	Consider the initial value problem
	\begin{equation} \label{eq:IVP}
		\begin{cases}
			u''(t) = f(t, u, u'),\ t \in [t_{0}, t_{1}] \\
			u(t_{0}) = u_{0},\ u'(t_{0}) = u'_{0}.
		\end{cases}
	\end{equation}
	Suppose that
	\begin{enumerate}[(1)]
		\item There exist a subsolution $u_{\mathrm{sub}}(t)$ and a supersolution
		      $u_{\mathrm{sup}}(t)$ such that
			  \begin{align*}
				 (i)\, &  u_{\mathrm{sub}}(t) \leqslant u_{\mathrm{sup}}(t),\ u_{\mathrm{sub}}'(t) \leqslant u_{\mathrm{sup}}'(t),\ t \in [t_{0}, t_{1}]; \\
		         (ii)\,	& u_{\mathrm{sub}}(t_{0}) \leqslant u_{\mathrm{sup}}(t_{0}),\ u_{\mathrm{sub}}'(t_{0}) \leqslant u_{\mathrm{sup}}'(t_{0}); \\
				 (iii)\, &\text{For any } t \in [t_{0}, t_{1}] )\text{ and } u \in [u_{\mathrm{sub}}(t), u_{\mathrm{sup}}(t)], \text{ we have} \\
				 &u_{\mathrm{sub}}''(t) \leqslant f(t, u, u_{\mathrm{sub}}'(t)),\ u_{\mathrm{sup}}''(t) \geqslant f(t, u, u_{\mathrm{sup}}'(t)).
			  \end{align*}
		\item The function $f$ on the right-hand side of the ODE in (\ref{eq:IVP}) 
		      is continuous on the domain $D = \{(t, u, v) \colon t_{0} \leqslant
              t \leqslant t_{1},\ u_{\mathrm{sub}}(t) \leqslant u \leqslant 
			  u_{\mathrm{sup}}(t),\ u_{\mathrm{sub}}'(t) \leqslant v \leqslant 
			  u_{\mathrm{sup}}'(t)\}$.
	\end{enumerate}
	Then the solution $u(t)$ of the initial value problem (\ref{eq:IVP}) satisfies
	$u(t) \in C^{2}[t_{0}, t_{1}]$ and $ u_{\mathrm{sub}}(t) \leqslant u(t)
	\leqslant u_{\mathrm{sup}}(t)$, $u_{\mathrm{sub}}'(t) \leqslant u'(t) \leqslant
	u_{\mathrm{sup}}'(t)$ for any $t \in [t_{0}, t_{1}]$.
\end{lemma}

To define subsolutions and  supersolutions of the ODE \eqref{eq:hill}, we introduce 
the upper and lower bound of $q(x_2)$ by 
\begin{equation*}
	q_{\max} \coloneqq \sup_{x_{2} \in [0, 2\pi]} q(x_{2}),\quad q_{\min} \coloneqq \inf_{x_{2} \in [0, 2\pi]} q(x_{2}),
\end{equation*}
respectively. Recalling the definition of $\eta_{n}$ in (\ref{eq:etan}), we estimate 
$\eta_{n}$ for large $\vert n\vert$ as follows.

\begin{lemma}\label{lem:estimateEtan}
	We have $\eta_{n} > 2$ if $\alpha_{n}^{2} > k^{2}q_{\max}$ and $\eta_{n} = 
	O(e^{2\pi \vert n\vert})$ as $\vert n\vert \to \infty$.
\end{lemma}

\begin{Proof}
    Using Lemma \ref{lem:compareOde}, we can estimate $\eta_{n}$ by constructing
	subsolutions and supersolutions for the ODEs (\ref{eq:hill}). The subsolution
	and supersolution with the initial conditions $w_{n,1}(0) = 1, w'_{n,1}(0) =
    0$ are given by
	\begin{align*}
		&u_{\mathrm{sub}}(x_{2}) = \frac{1}{2}\left(e^{\beta_{n}^{+}x_{2}} + e^{-\beta_{n}^{+}x_{2}}\right),\\
		&u'_{\mathrm{sub}}(x_{2}) = \frac{\beta_{n}^{+}}{2}\left(e^{\beta_{n}^{+}x_{2}} - e^{-\beta_{n}^{+}x_{2}}\right),\\
		&u_{\mathrm{sup}}(x_{2}) = \frac{1}{2}\left(e^{\beta_{n}^{-}x_{2}} + e^{-\beta_{n}^{-}x_{2}}\right),\\
		&u'_{\mathrm{sup}}(x_{2}) = \frac{\beta_{n}^{-}}{2}\left(e^{\beta_{n}^{-}x_{2}} - e^{-\beta_{n}^{-}x_{2}}\right),
	\end{align*}
	where $\beta_{n}^{+} = \sqrt{\alpha^{2}_{n} - k^{2}q_{\max}}$ and $\beta_{n}^{-}
    = \sqrt{\alpha^{2}_{n} - k^{2}q_{\min}}$. Similarly, the subsolution and supersolution
    with the initial condition $w_{n,2}(0) = 0, w'_{n,2}(0) = 1$ are
	\begin{align*}
		&u_{\mathrm{sub}}(x_{2}) = \frac{1}{2\beta_{n}^{+}}\left(e^{\beta_{n}^{+}x_{2}} - e^{-\beta_{n}^{+}x_{2}}\right),\\
		&u'_{\mathrm{sub}}(x_{2}) = \frac{1}{2}\left(e^{\beta_{n}^{+}x_{2}} + e^{-\beta_{n}^{+}x_{2}}\right),\\
		&u_{\mathrm{sup}}(x_{2}) = \frac{1}{2\beta_{n}^{-}}\left(e^{\beta_{n}^{-}x_{2}} - e^{-\beta_{n}^{-}x_{2}}\right),\\
		&u'_{\mathrm{sup}}(x_{2}) = \frac{1}{2}\left(e^{\beta_{n}^{-}x_{2}} + e^{-\beta_{n}^{-}x_{2}}\right).
	\end{align*}
    By Lemma \ref{lem:compareOde}, we have that
    \begin{align*}
        \frac{1}{2}\left(e^{\beta_{n}^{+}x_{2}} + e^{-\beta_{n}^{+}x_{2}}\right) &\leqslant w_{n,1}(x_{2}) \leqslant \frac{1}{2}\left(e^{\beta_{n}^{-}x_{2}} + e^{-\beta_{n}^{-}x_{2}}\right), \\
        \frac{1}{2}\left(e^{\beta_{n}^{+}x_{2}} + e^{-\beta_{n}^{+}x_{2}}\right) &\leqslant w'_{n,2}(x_{2}) \leqslant \frac{1}{2}\left(e^{\beta_{n}^{-}x_{2}} + e^{-\beta_{n}^{-}x_{2}}\right).
    \end{align*}
	
	For $\alpha_{n}^{2} > k^{2}q_{\max}$, the number $\beta_{n}^{+}$ is real-valued. Hence
	\begin{align*}
		\eta_{n} &= w_{n,1}(2\pi) + w'_{n,2}(2\pi) \\
		         &\geqslant e^{2\pi\beta_{n}^{+}} + e^{-2\pi\beta_{n}^{+}} \\
		         &> 2.
	\end{align*}
	Furthermore, it follows from Lemma \ref{lem:compareOde} with the above 
	subsolutions and supersolutions that $\eta_{n} = O(e^{2\pi \vert n\vert})$ 
	as $\vert n\vert \to +\infty$.
\end{Proof}

\begin{lemma} \label{lem:monoMu}
	Recall the definition of $\mu_{n, 1}$ in Definition \ref{def:characteristic}. 
	The sequences $\mu_{n, 1}$ and $\mu_{-n, 1}$ are monotonically increasing as 
	$\vert n\vert$ goes to $\infty$. Moreover, $\vert \mu_{n, 1}\vert = O(\vert n\vert)$.
\end{lemma}

\begin{Proof}
	We first prove that $\mu_{n, 1}$ is monotonically increasing as $n$ tends to 
	infinity. By (\ref{eq:polynomial}), it suffices to prove that
	$\eta_{n}$ is monotonically increasing as $n$ tends to infinity. We use the same
	argument as in the proof of Lemma \ref{lem:estimateEtan} to get the upper bound for $\eta_{n}$
	and the lower bound for $\eta_{n + 1}$, i.e., $\eta_{n} \leqslant e^{2\pi\beta_{n}^{-}}
	+ e^{-2\pi\beta_{n}^{-}}$ and $\eta_{n + 1} \geqslant e^{2\pi\beta_{n+1}^{+}}
    + e^{-2\pi\beta_{n+1}^{+}}$. If $e^{2\pi\beta_{n+1}^{+}}$ is greater than
    $e^{2\pi\beta_{n}^{-}}$, i.e., $n$ is greater than $\left(k^{2}(q_{\max}
    - q_{\min}) - 2 \alpha - 1\right)/2$, then $\eta_{n}$ is monotonically increasing.
    With the help of Definition \ref{def:characteristic} and (\ref{eq:polynomial}),
    we can show that $\mu_{n, 1}$ is monotonically increasing. Similarly, if $n$
    approaches the negative infinity, then we may have the upper bound for $\eta_{n}$
    and the lower bound for $\eta_{n - 1}$, i.e., $\eta_{n} \leqslant e^{2\pi\beta_{n}^{-}}
    + e^{-2\pi\beta_{n}^{-}}$ and $\eta_{n - 1} \geqslant e^{2\pi\beta_{n-1}^{+}}
    + e^{-2\pi\beta_{n-1}^{+}}$. Hence, if $n$ is less than $(k^{2}(q_{\min}
    - q_{\max}) - 2 \alpha + 1) / 2$, then $\mu_{n, 1}$ is  monotonically decreasing.
    By equation (\ref{eq:polynomial}) and $e^{2\pi \mu_{n, 1}} = \lambda_{n, 1}$,
    we know that $\vert \mu_{n, 1}\vert = O(\vert n\vert)$.
\end{Proof}

More precise estimates of $\eta_{n}$ and $\mu_{n,1}$ will be shown in the proof 
of Lemma \ref{lem:unibdd} below. These asymptotics enable us to justify the 
pointwise convergence of the proposed upward and downward radiation conditions.

\begin{theorem}[Pointwise Convergence]
	If $u \in H_{loc}^{1}(\mathbb{R}^{2})$, then the infinite series on the right-hand 
	side of \eqref{eq:urc} (resp.~\eqref{eq:drc}) in Definition \ref{def:radiationCondition} 
	is pointwisely convergent.
\end{theorem}

\begin{Proof}
    Without loss of generality, we only prove the theorem for the upward radiation
    condition. It suffices to prove that the radiation condition converges
    pointwisely in $x_{2} > d$ for $n$ sufficiently large.

    For any fixed $(x_{1}, x_{2})^{\mathrm{T}} \in \Omega_{d}^{+}$, there is a
    real number $\sigma \in (0, 2\pi]$ such that $v_{n, 2}(x_{2}) = v_{n, 2}(d
    + \sigma)$ for any $n \in \mathbb{Z}$. If $x_{2} > d + \sigma$, then we can
    decompose the upward mode into two parts. Namely, applying Young's inequality
    and Parseval's equality yields the estimate
    \begin{align*}
    	&2 \bigg\vert \sum_{\eta_{n} > 2} C_{n}^{+} e^{-\mu_{n, 1}x_{2}}v_{n, 2}(x_{2})e^{i \alpha_{n} x_{1}}\bigg\vert \\
    	\leqslant & \sum_{\eta_{n} > 2} 2 \bigg\vert C_{n}^{+} e^{-\mu_{n, 1}(d + \sigma)}v_{n, 2}(x_{2})e^{i \alpha_{n} x_{1}}e^{-\mu_{n, 1}(x_{2} - d - \sigma)}\bigg\vert \\
    	\leqslant & \sum_{\eta_{n} > 2} \bigg\vert C_{n}^{+} e^{-\mu_{n, 1}(d + \sigma)}v_{n, 2}(x_{2})e^{i \alpha_{n} x_{1}}\bigg\vert^{2} + \sum_{\eta_{n} > 2} \bigg\vert e^{-\mu_{n, 1}(x_{2} - d - \sigma)}\bigg\vert^{2} \\
    = & \sum_{\eta_{n} > 2} \bigg\vert C_{n}^{+} e^{-\mu_{n, 1}(d + \sigma)}v_{n, 2}(d+\sigma)e^{i \alpha_{n} x_{1}}\bigg\vert^{2}+ \sum_{\eta_{n} > 2} \bigg\vert e^{-\mu_{n, 1}(x_{2} - d - \sigma)}\bigg\vert^{2} \\
    	\leqslant & \,{\frac{1}{\tt p}}\int_{0}^{\tt p} \vert u(x_{1}, d+\sigma)\vert^{2} dx_{1} + \sum_{\eta_{n} > 2} \bigg\vert e^{-\mu_{n, 1}(x_{2} - d - \sigma)}\bigg\vert^{2} \\
    < & \;\infty,
    \end{align*}
    where the second term of the last inequality is convergent since $\mu_{n, 1}$
    is monotonically increasing and $\vert \mu_{n, 1}\vert = O(\vert n\vert)$
    (cf.~the proof of Lemma \ref{lem:monoMu}).
       
    If $x_{2} \in (d, d + \sigma]$, then we choose an appropriate $\epsilon$ 
	such that $d < \epsilon < x_{2}$ and $v_{n,2}(\epsilon) \neq 0$ (see Remark \ref{rem:nonzero}). 
	Using Parseval's equality and the asymptotic behavior of $v_{n,2}(x_{2})$ in
    Lemma \ref{lem:unibdd} yields that
    \begin{align*}
        &2 \bigg\vert \sum_{\eta_{n} > 2} C_{n}^{+} e^{-\mu_{n, 1}x_{2}}v_{n, 2}(x_{2})e^{i \alpha_{n} x_{1}}\bigg\vert \\
        \leqslant & \sum_{\eta_{n} > 2} 2 \bigg\vert C_{n}^{+} e^{-\mu_{n, 1}\epsilon}v_{n, 2}(x_{2})e^{i \alpha_{n} x_{1}}e^{-\mu_{n, 1}(x_{2} - \epsilon)}\bigg\vert \\
        \leqslant & \sum_{\eta_{n} > 2} \bigg\vert C_{n}^{+} e^{-\mu_{n, 1}\epsilon}v_{n, 2}(\epsilon)e^{i \alpha_{n} x_{1}}\bigg\vert^{2}\bigg\vert \frac{v_{n, 2}(x_{2})}{v_{n, 2}(\epsilon)}\bigg\vert^{2}+ \sum_{\eta_{n} > 2} \bigg\vert e^{-\mu_{n, 1}(x_{2} - \epsilon)}\bigg\vert^{2} .
    \end{align*}
    Indeed, in view of the subsequent (\ref{eq:vasymp}), we get
    $|\frac{v_{n, 2}(x_{2})}{v_{n, 2}(\epsilon)}|^2<C$ for $|n| > n_{a}$ with 
	sufficiently large $n_a$.
    \begin{align*}
        &2 \bigg\vert \sum_{|n|>n_a} C_{n}^{+} e^{-\mu_{n, 1}x_{2}}v_{n, 2}(x_{2})e^{i \alpha_{n} x_{1}}\bigg\vert \\
        \leqslant & \;\frac{C}{\tt p}\int_{0}^{\tt p} \vert u(x_{1}, \epsilon)\vert^{2} dx_{1} +\sum_{|n|>n_a} \bigg\vert e^{-\mu_{n, 1}(x_{2} - \epsilon)}\bigg\vert^{2}
        <  \;\infty.
    \end{align*}
    This completes the proof.
\end{Proof}


\section{Uniqueness and existence results} \label{sec:solve}

The aim of this section is to prove well-posedness of the diffraction problem
described in Section \ref{sec:problem} complemented with the classical upward
Rayleigh expansion in $x_2>d$ and the new downward radiation condition \eqref{eq:drc}
defined in Definition \ref{def:radiationCondition} in $x_{2} < b$. In Subsection
\ref{subsec:dtn} we carefully define the DtN maps on $\Gamma_{b}$ and $\Gamma_{d}$,
allowing us to truncate the physical domain to a bounded periodic cell. In Subsection
\ref{variation} we prove the Fredholm property of the resulting sesquilinear
form and then in Subsection \ref{well-posedness} the uniqueness and existence
of our diffraction problem. The solvability results depend on a well-defined
downward DtN mapping for the infinite inhomogeneous medium $q(x_{2})$ and the
mapping properties for its real and imaginary parts.

\subsection{Dirichlet-to-Neumann (DtN) maps} \label{subsec:dtn}

We first recall the classical $\hat{\alpha}$-quasiperiodic (upward) DtN map $T$ 
for $q \equiv 1$.
\begin{definition}[Classical upward DtN map]
    Given a quasi-periodic function $f$ over $\Gamma_{d}$, which has a Fourier 
	series expansion
    \begin{equation} \label{eq:FouSerExp}
        f (x_1)= \sum_{n \in \mathbb{Z}} f_{n} e^{i\alpha_{n} x_{1}},\quad
        f_{n} \coloneqq \frac{1}{\tt p} \int_{0}^{\tt p} f(x_{1}) e^{-i\alpha_{n}x_{1}} d x_{1}.
    \end{equation}
    The classical upward DtN map $T$ in a homogeneous medium (that is, $q(x_{2}) 
	\equiv 1$) is defined as
    \begin{equation}  \label{eq:stDtNmap}
        (Tf)(x_{1}) = \sum_{n \in \mathbb{Z}} i\beta_{n}f_{n}\;e^{i\alpha_{n}x_{1}},
    \end{equation}
    such that the relation $\partial_{x_{2}}u = Tu$ is equivalent with the classical 
    upward Rayleigh-expansion condition for $u$ in $x_{2} > d$. Note that this 
	classical upward Rayleigh-expansion condition is (\ref{eq:clurc}) with $-\beta_{n}$ 
	replaced by $+\beta_{n}$ if the domain is $\Omega_{b}^{-}$ instead of $\Omega_{d}^{+}$.
\end{definition}
Next we define the upward DtN map $T^{+}$ for the inhomogeneous periodic medium
$q(x_2)$. The downward DtN map $T^{-}$ can be treated similarly. In the remaining
part of this paper we make the following assumption on $q(x_2)$:  
\begin{description}
\item{\textbf{Assumption A:}} If there holds $\eta_{n} = 2$ for $\eta_{n}$ in 
(\ref{eq:etan}) and some $n \in \mathbb{Z}$, then the condition of Case (d) (ii) 
in Lemma \ref{lem:solvereal} is fulfilled. Similarly, if there holds  $\eta_{n} = -2$ 
for some $n \in \mathbb{Z}$, then the condition of Case (e) (ii) is fulfilled.
In other words, Case (d) (i) for $\eta_{n} = 2$ and Case (e) (i) for $\eta_{n} = -2$ 
are excluded.
\end{description}
Note that Assumption A is always satisfied if $q(x_2) \equiv 1$.

Consider the Dirichlet boundary value problem
\begin{equation} \label{eq:bvp2}
	\begin{cases}
		\Delta u(x_{1}, x_{2}) + k^{2}q(x_{2})u(x_{1}, x_{2}) = 0 &\ \text{ in } \Omega_{d}^{+},\\
		u = f \in H_{\hat{\alpha}}^{1/2}(\mathbb{R}) &\ \text{ on } \Gamma_{d},\\
		u \text{ satisfies the } \eqref{eq:urc} \text{ in Definition } \ref{def:radiationCondition}.
	\end{cases}
\end{equation}
The DtN map $T^{+}$ maps the Dirichlet value $f$ to the Neumann value of the solution 
$u$ of (\ref{eq:bvp2}), i.e.,
\begin{equation*}
	T^{+} \colon f \longmapsto \frac{\partial u(x_{1}, x_{2})}{\partial x_{2}}\bigg\vert_{\Gamma_{d}}.
\end{equation*}
Since $u$ satisfies the \eqref{eq:urc} and $f \in H_{\hat{\alpha}}^{1/2}(\mathbb{R})$, 
we have that (cf.~(\ref{eq:FouSerExp}))
\begin{equation*}
	u(x_{1}, d) = \sum_{n \in \mathbb{Z}} C_{n}^{+}u_{n}^{+}(d) e^{i\alpha_{n}x_{1}} = f(x_{1}) = \sum_{n \in \mathbb{Z}} f_{n} e^{i\alpha_{n}x_{1}}\quad \text{ for all } x_{1} \in \mathbb{R}.
\end{equation*}
If $u_{n}^{+}(d) \neq 0$ for all $n$, then we can compute $C_{n}^{+} = f_{n}/u_{n}^{+}(d)$
and formally rewrite the unique solution $u(x_{1}, x_{2})$ as
\begin{equation} \label{eq:rayleigh}
	u(x_{1}, x_{2}) = \sum_{n \in \mathbb{Z}} \frac{f_{n}}{u_{n}^{+}(d)} u_{n}^{+}(x_{2})e^{i\alpha_{n}x_{1}}.
\end{equation}
Taking the derivative of (\ref{eq:rayleigh}) w.r.t. $x_{2}$ on $\Gamma_{d}$, we get
\begin{equation*} 
	\frac{\partial u(x_{1}, x_{2})}{\partial x_{2}}\bigg\vert_{\Gamma_{d}} = \sum_{n \in \mathbb{Z}} \frac{{u_{n}^{+}}'(d)}{u_{n}^{+}(d)} f_{n}\; e^{i\alpha_{n}x_{1}}.
\end{equation*}

\begin{remark} \label{rem:nonzero}
\begin{itemize}
	\item[(i)] It is possible to choose $d \in \mathbb{R}$ such that $u_{n}^{+}(d) \neq 0$ 
	for all $n\in \mathbb{N}$. For this purpose we need to consider the zeros of 
	solutions of equation (\ref{eq:hill}). Property 4.1.2 in \cite[Section 4.1]{Ea1} 
	asserts that any solution of (\ref{eq:hill}) has a finite number of zeros only 
	in any bounded closed interval. Furthermore, Theorem 4.23 in \cite[Section 4.2]{Ea1} 
	shows that if $n$ is sufficiently large such that $k^{2}q(x_{2}) - \alpha_{n}^{2} \leqslant 0$, 
	then no solutions of (\ref{eq:hill}) has more than one zero. Let $\mathcal{Z}$ 
	denote the set of all zeros of solutions to (\ref{eq:hill}) for all $n \in \mathbb{N}$. 
	Then $\mathcal{Z}$ is a countable set. Hence, we can always find an appropriate 
	$d \in \mathbb{R} \backslash \mathcal{Z}$ such that $u_{n}^{+}(d) \neq 0$ for 
	all $n \in \mathbb{Z}$. Below we shall assume that $u_{n}^{+}(d) \neq 0$ for 
	all $n \in \mathbb{Z}$. If $u_{n}^{+}(d)=0$ but if there exists the finite limit
	$ \lim_{\,0<\varepsilon\rightarrow 0} \frac{{u_{n}^{+}}'(d+\varepsilon)}{u_{n}^{+}(d+\varepsilon)}$, 
	then $\frac{{u_{n}^{+}}'(d)}{u_{n}^{+}(d)}$ in (\ref{eq:newDtN}) can be replaced 
	by $ \lim_{\,0<\varepsilon\rightarrow 0} \frac{{u_{n}^{+}}'(d+\varepsilon)}{u_{n}^{+}(d+\varepsilon)}$.

	\item[(ii)] If case (d) (i) or case (e) (i) occurs, then the DtN map $T^{+}$ 
	is not well-defined, because solutions to the Dirichlet problem \eqref{eq:bvp2} 
	are not unique. For example, if $f(x_1)=f_{n} e^{i\alpha_{n} x_{1}}$ with 
	$\eta_{n} = 2$ for some $n \in \mathbb{N}$, then any function of the form 
	$u(x)= (c_{1} v_{n,1}(x_{2}) + c_{2} v_{n,2}(x_{2})) e^{i\alpha_{n} x_{1}}$, 
	with coefficients $c_j\in \mathbb{C}$ satisfying $c_1 v_{n,1}(d)+c_2 v_{n,2}(d)=f_n$, 
	is a solution of the Dirichlet problem \eqref{eq:bvp2}.
	
	\item[(iii)] The assumption $u_{n}^{+}(d) \neq 0$ implies that $v_{n, 2}(d) \neq 0$ 
	for the case $\vert\eta_n\vert>2$ (cf.~Definition \ref{def:modes}). If $q(x_{2}) 
	\equiv 1$, then we have $u_{n}^{+}(x_{2})= e^{i\beta_{n} x_{2}}$ and $u_{n}^{+}(d) 
	\neq 0$ for all $n$ and $d$.
\end{itemize}
\end{remark}

The upward and downward DtN maps $T^\pm$ for an inhomogeneous medium $q(x_2)$ 
are defined as follows.
\begin{definition}[DtN maps for $q(x_2)$] \label{def:dtnMapping}
	Choose $d\in \mathbb{R}$ such that $u_{n}^{+}(d) \neq 0$ for all $n \in \mathbb{Z}$. 
	For quasi-periodic functions $f \in H_{\hat{\alpha}}^{1/2}(\Gamma_{d})$, we 
	define the DtN map $T^{+} \colon H_{\hat{\alpha}}^{1/2}(\Gamma_{d}) \to H_{\hat{\alpha}}^{-1/2}(\Gamma_{d})$ 
	by
	\begin{equation} \label{eq:newDtN}
		(T^{+} f)(x_{1}) \coloneqq \sum_{n \in \mathbb{Z}} \frac{{u_{n}^{+}}'(d)}{u_{n}^{+}(d)} f_{n}\; e^{i\alpha_{n} x_{1}},
	\end{equation}
	where $f\vert_{\Gamma_{d}}$ has the Fourier series expansion (\ref{eq:FouSerExp}).
	
    Choosing $b\in \mathbb{R}$ such that $u_{n}^{-}(b) \neq 0$ for all $n \in \mathbb{Z}$,
    we define the downward Dirichlet-to-Neumann map $T^{-}$ on $\Gamma_{b}$ by
	\begin{equation} \label{eq:newDtNminus}
		(T^{-}f)(x_{1}) \coloneqq \sum_{n \in \mathbb{Z}} -\frac{{u_{n}^{-}}'(b)}{u_{n}^{-}(b)} f_{n}\; e^{i\alpha_{n}x_{1}},
	\end{equation}
	where $f_{n} \in \mathbb{C}$ denotes the Fourier coefficients of $f$ on $\Gamma_{b}$.
\end{definition}

To get the boundedness of the DtN maps $T^{\pm}$ we need to investigate the 
asymptotic behavior of $v_{n,j}$ ($j = 1, 2$) as $|n| \rightarrow \infty$. Below 
we shall only focus on the upward DtN mapping $T^{+}$. By Lemma \ref{lem:estimateEtan}, 
we have $\eta_{n} > 2$ if $|n|$ is sufficiently large.

\begin{lemma}[Boundedness] \label{lem:unibdd}
	Let $v_{n,2}$ be defined as in Lemma \ref{lem:solvereal}. There exists $M
    > 0$ such that $\big\vert \frac{v_{n,2}'(d)}{v_{n,2}(d)}\big\vert < M$ for
    all $\eta_{n} > 2$.
\end{lemma}

\begin{Proof}
\noindent The proof consists of three steps:
\begin{enumerate} [1.]
    \item Get the solutions of (\ref{eq:hill}) with the initial conditions (\ref{eq:initial}).
    \item Derive the asymptotic behavior of the characteristic exponents defined in Definition \ref{def:characteristic}.
    \item Derive the asymptotic behavior of $v_{n, j}$, $j = 1, 2$, according to the modes defined in Definition \ref{def:modes}.
\end{enumerate}
For simplicity, we abuse notation by dropping the subscript $n$ of $\alpha_{n}$, 
$\eta_n$, $u_n$ and $\lambda_{n,j}$ in the following steps.

\textbf{Step 1: Fundamental solutions.} Let $\psi$ be the solutions of $u''(t) -
\alpha^{2}u(t) = 0$ which satisfy the initial conditions
\begin{equation} \label{eq:initial1}
	u(0) = 1,\ u'(0) = 0,
\end{equation}
or
\begin{equation} \label{eq:initial2}
	u(0) = 0,\ u'(0) = 1.
\end{equation}
Hence we have that
\begin{equation*}
	\psi(t) =
	\begin{cases}
		\frac{1}{2}(e^{\alpha t} + e^{-\alpha t}), &\text{ with condition } (\ref{eq:initial1}), \\
		\frac{1}{2\alpha}(e^{\alpha t} - e^{-\alpha t}), &\text{ with condition } (\ref{eq:initial2}).
	\end{cases}
\end{equation*}
We look for solutions to (\ref{eq:hill}) of the form $u(t) =
\psi(t)\varphi(t)$. Therefore, $\varphi(t)$ should satisfy the
equation
\begin{equation*}
	\varphi''(t) + H(t)\varphi'(t) + k^{2}q(t)\varphi(t) = 0,
\end{equation*}
where
\begin{equation*}
	H(t) \coloneqq 2\frac{\psi'(t)}{\psi(t)} =
	\begin{cases}
		2\alpha\tanh(\alpha t)&\text{ with condition } (\ref{eq:initial1}), \\
		2\alpha\cosh(\alpha t)&\text{ with condition } (\ref{eq:initial2}).
	\end{cases}
\end{equation*}
The initial conditions (\ref{eq:initial1}) and (\ref{eq:initial2}) for $u$ are satisfied if and only if $\varphi$ satisfies the initial conditions $\varphi(0) = 1$ and $\varphi'(0) = 0$.
We get $H(t) = 2V'(t)/V(t)$, where
\begin{equation*}
	V(t) \coloneqq
	\begin{cases}
		\cosh(\alpha t)&\text{ with condition } (\ref{eq:initial1}), \\
		\sinh(\alpha t)&\text{ with condition } (\ref{eq:initial2}).
	\end{cases}
\end{equation*}
By straightforward computation, we have that
\begin{equation*}
    [V^{2}(t)\varphi'(t)]' = V^{2}(t)
(\varphi''(t) + H(t) \varphi'(t)) = -V^{2}(t)k^{2}q(t)\varphi(t).
\end{equation*}
Integrating both sides w.r.t.~$t$ from $0$ to $t$, we get
\begin{equation*}
    \varphi'(t) = -k^{2}V^{-2}(t) \int_{0}^{t} V^{2}(\tau)
q(\tau)\varphi(\tau)\dif\tau.
\end{equation*}
Integrating $\varphi'(t)$ again w.r.t.~$t$ from $0$ to $t$, we obtain
\begin{equation} \label{eq:varphi}
	\begin{aligned}
		\varphi(t) &= \varphi(0) + \int_{0}^{t} \varphi'(\tau)\dif\tau \\
		           &= 1 - k^{2} \int_{0}^{t} V^{-2}(\tau)\int_{0}^{\tau} V^{2}(\sigma)q(\sigma)\varphi(\sigma)\dif\sigma \dif\tau \\
		           &= 1 - k^{2} \int_{0}^{t} q(\sigma)V^{2}(\sigma)\left(\int_{\sigma}^{t} V^{-2}(\tau)\dif\tau\right) \varphi(\sigma)\dif\sigma \\
		           &= 1 - k^{2} \int_{0}^{t} \Big[q(\sigma)V^{2}(\sigma)(\tilde{H}(t) - \tilde{H}(\sigma))\Big] \varphi(\sigma)\dif\sigma,
	\end{aligned}
\end{equation}
where
\begin{equation*}
	\tilde{H}(t) \coloneqq \frac{1}{\alpha}
	\begin{cases}
		\tanh(\alpha t)&\text{ with condition } (\ref{eq:initial1}), \\
		\coth(\alpha t)&\text{ with condition } (\ref{eq:initial2}).
	\end{cases}
\end{equation*}

From the discussions above, we have an integral operator $K$ with kernel $k(t,s)$
defined by
\begin{equation*}
	(K\varphi)(t) = \int_{0}^{t} k(t,\sigma)\varphi(
	\sigma)\dif\sigma,
\end{equation*}
where $k(t,\sigma) = q(\sigma)V^{2}(\sigma)(\tilde{H}(t) - \tilde{H}(\sigma))$.
We note that the expression of the kernel $k(t,\sigma)$ relies on the initial
conditions \eqref{eq:initial1} and \eqref{eq:initial2}. Now we can formulate
(\ref{eq:varphi}) as the operator equation
\begin{equation} \label{eq:opequ}
	\varphi(t) + k^{2}K\varphi(t) = 1.
\end{equation}
An expansion of the solution $\varphi$ can be established by the Neumann series.
Firstly, We need to estimate the norm of the operator $K$ restricted to a finite
interval $[0, T]$, $T > 2\pi$.

We observe that $\vert k(t, \sigma) \vert \leqslant \vert q(\sigma)\vert / \alpha$.
Indeed, with the condition (\ref{eq:initial1}), we have that
\begin{align*}
	\vert k(t,\sigma)\vert &\leqslant \frac{1}{4\alpha} \vert q(\sigma)\vert (e^{\alpha\sigma} + e^{-\alpha\sigma})^{2} \left(1 - \frac{e^{\alpha\sigma} - e^{-\alpha\sigma}}{e^{\alpha\sigma} + e^{-\alpha\sigma}}\right)\\
	                       &\leqslant \frac{1}{2\alpha} \vert q(\sigma)\vert (e^{\alpha\sigma} + e^{-\alpha\sigma})e^{-\alpha\sigma}\\
	                       &\leqslant \frac{1}{\alpha} \vert q(\sigma)\vert.
\end{align*}
With the condition (\ref{eq:initial2}), we have that
\begin{align*}
	\vert k(t,\sigma)\vert &\leqslant \frac{1}{4\alpha} \vert q(\sigma)\vert (e^{\alpha\sigma} - e^{-\alpha\sigma})^{2} {\left(\frac{e^{\alpha\sigma} + e^{-\alpha\sigma}}{e^{\alpha\sigma} - e^{-\alpha\sigma}} - 1\right)} \\
	                       &\leqslant \frac{1}{2\alpha} \vert q(\sigma)\vert (e^{\alpha\sigma} - e^{-\alpha\sigma}) e^{-\alpha\sigma} \leqslant \frac{1}{2\alpha} \vert q(\sigma)\vert, \\
	\vert \partial_tk(t,\sigma)\vert &{\leqslant \frac{1}{4\alpha} \vert q(\sigma)\vert (e^{\alpha\sigma} - e^{-\alpha\sigma})^{2} \left|\partial_t\left(\frac{e^{\alpha t} + e^{-\alpha t}}{e^{\alpha t} - e^{-\alpha t}}\right)\right|} \\
	                                 &{\leqslant \vert q(\sigma)\vert \left|\frac{e^{\alpha\sigma} - e^{-\alpha\sigma}}{e^{\alpha t} - e^{-\alpha t}}\right|^2\leqslant \vert q(\sigma)\vert}
\end{align*}
for $0 < s \leqslant t < T$. Hence,
\begin{equation*}
  \Vert K\Vert \leqslant C \sup \vert k(t,\sigma)\vert \leqslant C \frac{1}{\alpha}
\end{equation*}
{and, for the case (\ref{eq:initial2}), we get
$\Vert \partial_t K\Vert\leqslant C$.}

For sufficiently large $\alpha$, $(I + k^{2}K)^{-1}$ exists and then
\begin{equation} \label{eq:expanPhi}
	\varphi(t) = \sum_{m = 0}^{\infty}\big((-k^{2}K)^{m} 1 \big)(t) = 1 - k^{2}(K1)(t) + O\!\left(\frac{1}{\alpha^{2}}\right),\quad 0\leqslant t \leqslant T.
\end{equation}
Similarly, in case of (\ref{eq:initial2}), we get
\begin{equation*}
  \varphi'(t) = - k^{2}(\partial_tK1)(t) + k^{4}\partial_tK(K1)(t) + O\!\left(\frac{1}{\alpha^{2}}\right),\quad 0 \leqslant t \leqslant T.
\end{equation*}

\textbf{Step 2: Characteristic exponents.} Using the fundamental solutions in 
Step 1 and (\ref{eq:etan}), we have that (ignoring higher order terms)
\begin{align*}
    \eta &= \frac{1}{2}(e^{2\pi\alpha} + e^{-2\pi\alpha})\left(1 - O\!\left(\frac{1}{\alpha}\right)\right) + \frac{1}{2}(e^{2\pi\alpha} + e^{-2\pi\alpha})\left(1 - O\!\left(\frac{1}{\alpha}\right)\right)  +\;O\!\left(\frac{1}{\alpha}\right)\\
         &= e^{2\pi\alpha}\left(1 - O\!\left(\frac{1}{\alpha}\right)\right).
\end{align*}

Then by (\ref{eq:polynomial}) we can get
\begin{align*}
    \lambda_{1} &= \frac{\eta}{2} + \sqrt{\frac{\eta^{2}}{4} - 1} = e^{2\pi\alpha}\left(1 - O\left(\frac{1}{\alpha}\right)\right), \\
    \lambda_{2} &= \frac{1}{\lambda_{1}} = e^{-2\pi\alpha}\left(1 + O\left(\frac{1}{\alpha}\right)\right).
\end{align*}

It follows from Definition \ref{def:characteristic} that the asymptotic behavior of
characteristic exponents is
\begin{equation*}
    \mu_{1} = \frac{1}{2\pi} \log \lambda_{1}
            = \alpha + O\left(\frac{1}{\alpha}\right) \text{ and } \mu_{2} = -\mu_{1} = -\alpha + O\left(\frac{1}{\alpha}\right).
\end{equation*}
Furthermore, we know that
\begin{equation*}
    \mu_{j}^{2} - \alpha^{2} = (\mu_{j} - \alpha)(\mu_{j} + \alpha) = O\left(\frac{1}{\alpha}\right) O\left(2\alpha + O\left(\frac{1}{\alpha}\right)\right) = O(1),\ \text{for } j = 1, 2.
\end{equation*}

\textbf{Step 3: Asymptotics of Hill's equation.}  According to the case (a) in
Definition \ref{def:modes}, we assume that $u(x_{2}) = e^{-\mu x_{2}}v(x_{2})$, where
$v(0) = v(2\pi)$ and $v'(0) = v'(2\pi)$. Substituting this representation of $u$ into equation
(\ref{eq:hill}), we obtain
\begin{equation*}
    e^{-\mu t}\left(v''(t) - 2\mu v'(t) + (k^{2}q(t) + \mu^{2} - \alpha^{2})v(t)\right) = 0.
\end{equation*}
Then we have that $(e^{-2\mu t}v'(t))' = -e^{-2\mu t}(k^{2}q(t)+\mu^{2}-\alpha^{2})v(t)$. Hence
\begin{align}
  e^{-2\mu t} v'(t) &= v'(0) - \int_{0}^{t} e^{-2\mu\tau} (k^{2}q(\tau)+\mu^{2}-\alpha^{2})v(\tau)\dif\tau, \label{eq:eq1} \\
  e^{-4\pi\mu } v'(2\pi) &=  e^{-2\mu t}v'(t) - \int_{t}^{2\pi} e^{-2\mu\tau} (k^{2}q(\tau)+\mu^{2}-\alpha^{2})v(\tau)\dif\tau. \label{eq:eq2}
\end{align}
Multiplying (\ref{eq:eq1}) by $e^{-4\mu\pi}$ and substituting $v'(2\pi)=v'(0)$ in (\ref{eq:eq2}), respectively, we get the two identities
\begin{equation} \label{eq:derivativeV1}
    e^{-4\pi\mu}e^{-2\mu t} v'(t) = e^{-4\pi\mu}v'(0) - e^{-4\pi\mu}\int_{0}^{t} e^{-2\mu\tau} (k^{2}q(\tau)+\mu^{2}-\alpha^{2})v(\tau)\dif\tau
\end{equation}
and
\begin{equation} \label{eq:derivativeV2}
    e^{-2\mu t} v'(t) = e^{-4\pi\mu} v'(0) + \int_{t}^{2\pi}e^{-2\mu\tau} [k^{2}q(\tau)+\mu^{2}-\alpha^{2}]v(\tau)\dif\tau.
\end{equation}
Subtracting (\ref{eq:derivativeV1}) from (\ref{eq:derivativeV2}) gives
\begin{align*}
    v'(t) = &\ \frac{1}{1-e^{-4\pi\mu}} \int_{t}^{2\pi} e^{2\mu(t-\tau)}(k^{2}q(\tau)+\mu^{2}-\alpha^{2})v(\tau)\dif\tau \\
            & +\frac{e^{-4\pi\mu}}{1-e^{-4\pi\mu}} \int_{0}^{t} e^{2\mu(t-\tau)}(k^{2}q(\tau)+\mu^{2}-\alpha^{2})v(\tau)\dif\tau \\
          = &\ \frac{e^{4\pi\mu}}{e^{4\pi\mu}-1} \int_{t}^{2\pi} e^{2\mu(t-\tau)}(k^{2}q(\tau)+\mu^{2}-\alpha^{2})v(\tau)\dif\tau \\
            & +\frac{1}{e^{4\pi\mu}-1} \int_{0}^{t} e^{2\mu(t-\tau)}(k^{2}q(\tau)+\mu^{2}-\alpha^{2})v(\tau)\dif\tau.
\end{align*}

There are two cases to consider. We first suppose that $v(0) \neq 0$. Without loss
of generality, we suppose $v(0) = 1$. Then
\begin{align*}
    v(t) = &\ 1 + \int_{0}^{t} v'(\sigma) \dif\sigma \\
         = &\ 1 + \frac{e^{4\pi\mu}}{e^{4\pi\mu}-1} \int_{0}^{t}\int_{\sigma}^{2\pi} e^{2\mu(\sigma-\tau)}(k^{2}q(\tau)+\mu^{2}-\alpha^{2})v(\tau) \dif\tau \dif\sigma \\
           & + \frac{1}{e^{4\pi\mu}-1}\int_{0}^{t}\int_{0}^{\sigma} e^{2\mu(\sigma-\tau)}(k^{2}q(\tau)+\mu^{2}-\alpha^{2})v(\tau) \dif\tau \dif\sigma \\
         = &\ 1 + \frac{e^{4\pi\mu}}{e^{4\pi\mu}-1} \int_{0}^{2\pi} \int_{0}^{\min(t,\tau)} e^{2\mu(\sigma-\tau)} \dif\sigma (k^{2}q(\tau)+\mu^{2}-\alpha^{2})v(\tau)\dif\tau \\
           & + \frac{1}{e^{4\pi\mu}-1} \int_{0}^{t} \int_{\tau}^{t} e^{2\mu(\sigma-\tau)} \dif\sigma (k^{2}q(\tau)+\mu^{2}-\alpha^{2})v(\tau) \dif\tau \\
         = &\ 1 + \frac{e^{4\pi\mu}}{e^{4\pi\mu}-1} \int_{0}^{2\pi} \frac{1}{2\mu}(e^{2\mu\min(t-\tau,0)} - e^{-2\mu\tau}) (k^{2}q(\tau)+\mu^{2}-\alpha^{2})v(\tau) \dif\tau \\
           & + \frac{1}{e^{4\pi\mu}-1} \int_{0}^{t} \frac{1}{2\mu}(e^{2\mu(t-\tau)} - 1) (k^{2}q(\tau)+\mu^{2}-\alpha^{2})v(\tau) \dif\tau \\
         = &\ 1 + \frac{1}{2\mu} \int_{0}^{2\pi}\left(e^{2\mu\min(t-\tau,0)} + \frac{e^{2\mu(t-\tau)} - e^{2\mu(2\pi-\tau)}}{e^{4\pi\mu}-1}\right) (k^{2}q(\tau)+\mu^{2}-\alpha^{2})v(\tau) \dif\tau.
\end{align*}
We define an integral operator $L$ by
\begin{equation*}
    (Lv)(t) \coloneqq \int_{0}^{2\pi} l(t,\tau) v(\tau) \dif\tau,
\end{equation*}
where the kernel $l(t,\tau)$ is defined by
\begin{equation*}
    l(t, \tau) = \frac{1}{2} \left(e^{2\mu\min(t-\tau,0)} + \frac{e^{2\mu(t-\tau)} - e^{2\mu(2\pi-\tau)}}{e^{4\pi\mu}-1}\right) (k^{2}q(\tau)+\mu^{2}-\alpha^{2}).
\end{equation*}
Then we obtain the integral equation
\begin{equation*}
    v(t) - \frac{1}{\mu}(Lv)(t) = 1.
\end{equation*}
We observe that $\vert l(t,\tau)\vert \leqslant C$. By the Neumann series, we can show that
\begin{equation}  \label{eq:vasymp}
    v(t) = \sum_{m = 0}^{\infty} \left(\frac{1}{\mu}L\right)^{m} 1 = 1 + \frac{1}{\mu}\int_{0}^{2\pi} l(t,\tau) \dif\tau + O\left(\frac{1}{\mu^{2}}\right).
\end{equation}

In the second case, if $v(0) = 0$, then the integral equation becomes $v(t) - 
\frac{1}{\mu}(Lv)(t) = 0$. Taking the norm of both sides of the equation 
$\mu v(t) = (Lv)(t)$, we find that $\mu \leqslant \Vert L\Vert \leqslant C$.
This implies that there are only finitely many $\mu = \mu_{n, 1}$ with 
$v(0) = 0$ (cf. Lemma \ref{lem:monoMu}).

It remains to estimate the derivative $v'(t)$. We first compute $\partial_{t}\frac{1}{\mu}(Lv)(t)$ by
\begin{align*}
    \partial_{t}\frac{1}{\mu}(Lv)(t) = & \ \partial_{t} \frac{1}{\mu}\int_{0}^{2\pi} l(t,\tau) v(\tau) \dif\tau \\
                                     = & \ \partial_{t} \frac{1}{2\mu} \int_{0}^{t} (k^{2}q(\tau)+\mu^{2}-\alpha^{2})v(\tau) \dif\tau \\
                                       & + \partial_{t} \frac{1}{2\mu} \int_{t}^{2\pi} e^{2\mu(t-\tau)}(k^{2}q(\tau)+\mu^{2}-\alpha^{2})v(\tau) \dif\tau \\
                                       & + \partial_{t} \frac{1}{2\mu} \int_{0}^{2\pi}  \frac{e^{2\mu(t-\tau)} - e^{2\mu(2\pi-\tau)}}{e^{4\pi\mu}-1} (k^{2}q(\tau)+\mu^{2}-\alpha^{2})v(\tau)  \dif\tau \\
                                     = & \ \int_{t}^{2\pi} e^{2\mu(t-\tau)}(k^{2}q(\tau)+\mu^{2}-\alpha^{2})v(\tau) \dif\tau \\
                                       & + \int_{0}^{2\pi} \frac{e^{2\mu(t-\tau)}}{e^{4\pi\mu}-1} (k^{2}q(\tau)+\mu^{2}-\alpha^{2})v(\tau) \dif\tau.
\end{align*}
It follows that
\begin{align*}
    \frac{\partial_{t}v(t)}{v(t)} = \frac{\partial_{t}\frac{1}{\mu}(Lv)(t)}{v(t)}
                                  = & \ \frac{1}{v(t)}\int_{t}^{2\pi} e^{2\mu(t-\tau)}(k^{2}q(\tau)+\mu^{2}-\alpha^{2})v(\tau) \dif\tau \\
                                    & + \frac{1}{v(t)}\int_{0}^{2\pi} \frac{e^{2\mu(t-\tau)}}{e^{4\pi\mu}-1} (k^{2}q(\tau)+\mu^{2}-\alpha^{2})v(\tau) \dif\tau.
\end{align*}
Hence, for sufficiently large $\mu$ (cf.~(\ref{eq:vasymp})),
\begin{equation*}
    \bigg\Vert \frac{\partial_{t}v(t)}{v(t)}\bigg\Vert_{L^{\infty}} \leqslant C.
\end{equation*}

Recalling the upward modes $u_{n}^{+}(x_{2}) = e^{-\mu_{n, 1}x_{2}}v_{n, 2}(x_{2})$ 
for $\eta_{n} > 2$ defined in Definition \ref{def:modes}, we need to estimate 
$v_{n,2}$ and $v_{n,2}'$. Similarly to the discussions in Step 2, we know that 
$\mu_{n,1}^{2} - \alpha_{n}^{2} = O(1)$ as $n$ tends to infinity. Since $v_{n,2}$ 
is a periodic function, we have $v_{n,2}(0) = v_{n,2}(2\pi)$, $v_{n,2}'(0) = v_{n,2}'(2\pi)$. 
Thus we can apply the methods in Step 3 to estimate $v_{n,2}$ and $v_{n,2}'$. 
Hence we have that for sufficiently large $n$ there exists a constant M such that
\begin{equation*}
    \bigg\vert \frac{v_{n,2}'(d)}{v_{n,2}(d)}\bigg\vert \leqslant \bigg\Vert \frac{v_{n,2}'(t)}{v_{n,2}(t)}\bigg\Vert_{L^{\infty}} \leqslant M.
\end{equation*}
\end{Proof}

\begin{theorem} \label{thm:continuousDtN}
	The DtN map $T^{+} \colon H_{\hat{\alpha}}^{1/2}(\mathbb{R}) \to
	H_{\hat{\alpha}}^{-1/2}(\mathbb{R})$ is continuous, i.e., there exists
	a positive constant $C$ such that
	\begin{equation} \label{eq:dtnbdd}
		\Vert T^{+}f\Vert_{H_{\hat{\alpha}}^{-1/2}(\mathbb{R})} \leqslant C\Vert f\Vert_{H_{\hat{\alpha}}^{1/2}(\mathbb{R})} \quad \text{for all } f \in H_{\hat{\alpha}}^{1/2}(\mathbb{R}).
	\end{equation}
\end{theorem}

\begin{Proof}
	For any $f \in H_{\hat{\alpha}}^{1/2}(\mathbb{R})$, we obtain that
	\begin{align*}
		\Vert T^{+} f\Vert_{H_{\hat{\alpha}}^{-1/2}(\mathbb{R})}^{2} &= {\tt p} \sum_{n \in \mathbb{Z}} (1 + \alpha_{n}^{2})^{-1/2} \bigg\vert \frac{{u_{n}^{+}}'(d)}{u_{n}^{+}(d)}\bigg\vert^{2} \vert f_{n}\vert^{2} \\
		        &= {\tt p} \sum_{n \in \mathbb{Z}} (1 + \alpha_{n}^{2})^{1/2} \bigg[(1 + \alpha_{n}^{2})^{-1/2}\bigg\vert\frac{{u_{n}^{+}}'(d)}{u_{n}^{+}(d)}\bigg\vert\bigg]^{2} \vert f_{n}\vert^{2}.
	\end{align*}
	By Lemma \refeq{lem:estimateEtan}, it suffices to prove that (\ref{eq:dtnbdd})
	holds for $T^{+}$ with a summation running over all $n$ with $\eta_{n} > 2$, i.e., 
	with $\vert n\vert$ sufficiently large. If follows from $\vert \mu_{n, 1}\vert 
	= O(\vert n\vert)$ (cf.~Lemma \ref{lem:monoMu}) and the boundedness of 
	$\frac{v_{n,2}'(d)}{v_{n,2}(d)}$ (cf.~Lemma \ref{lem:unibdd}) that
	\begin{equation*}
		(1 + \alpha_{n}^{2})^{-1/2}\bigg\vert\frac{{u_{n}^{+}}'(d)}{u_{n}^{+}(d)}\bigg\vert = \frac{\bigg\vert -\mu_{n, 1} + \frac{v_{n,2}'(d)}{v_{n,2}(d)}\bigg\vert}{(1 + \alpha_{n}^{2})^{1/2}} \leqslant C \text{ for all } \eta_{n} > 2.
	\end{equation*}
	Hence, we have that for any $f \in H_{\hat{\alpha}}^{1/2}\left(\mathbb{R}\right)$
	there exists a constant $C > 0$ such that
	\begin{equation*}
		\Vert T^{+}f\Vert_{H_{\hat{\alpha}}^{-1/2}(\mathbb{R})} \leqslant C\Vert f\Vert_{H_{\hat{\alpha}}^{1/2}(\mathbb{R})},
	\end{equation*}
	which completes the proof.
\end{Proof}

\begin{remark}
	Lemma \ref{lem:unibdd} and Theorem \ref{thm:continuousDtN} also hold for the downward DtN map $T^{-}$.
\end{remark}

\subsection{Variational formulation} \label{variation}

\begin{figure}[htp]
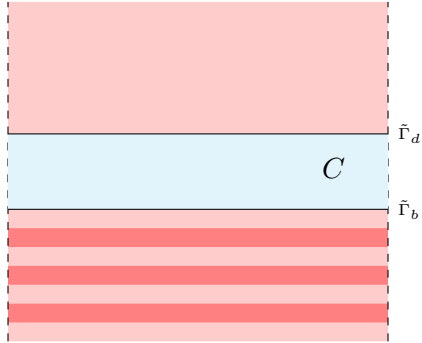

	\centering
	\includestandalone{figs/bvpc}
	\caption{Illustration of a periodic cell $C \coloneqq (0, {\tt p}) \times (b,d)$, which is covered by a homogeneous medium in $x_2>d$ and sits above an inhomogeneous half plane $x_2<b$. }
	\label{fig:periodic}
\end{figure}

To truncate the computational domain for the scattering problem to a bounded periodic 
cell (see Figure \ref{fig:periodic}), we introduce the bounded domain
\begin{equation*}
    C \coloneqq \{(x_{1}, x_{2})^{\mathrm{T}} \in \mathbb{R}^{2} : 0 < x_{1} < {\tt p},\ b < x_{2} < d\},
\end{equation*}
and the boundaries
\begin{align*}
    \tilde{\Gamma}_{d} \coloneqq \{(x_{1}, x_{2})^{\mathrm{T}} \in \mathbb{R}^{2} : 0 < x_{1} < {\tt p},\ x_{2} = d\}, \\
    \tilde{\Gamma}_{b} \coloneqq \{(x_{1}, x_{2})^{\mathrm{T}} \in \mathbb{R}^{2} : 0 < x_{1} < {\tt p},\ x_{2} = b\}.
\end{align*}

We consider the downward DtN map $T^{-}$ (cf.~Definition \ref{def:dtnMapping})
on $\tilde{\Gamma}_{b}$ and the standard upward DtN map $T$ on $\tilde{\Gamma}_{d}$ 
(cf.~(\ref{eq:stDtNmap})). Formulating the radiation conditions over $\tilde{\Gamma}_{d}$ 
and $\tilde{\Gamma}_{b}$ with the help of the DtN maps, we consider the boundary value problem
\begin{equation*}
	\begin{cases}
        \Delta u(x_{1}, x_{2}) + k^{2}\tilde{q}(x_{1},x_{2}) u(x_{1}, x_{2})
        = 0 &\;(x_1,x_2)^{\mathrm{T}}\in  C, \\[0.2cm]
        {\left[Tu{(\cdot,d)}\right](x_1)} \;\;=\phantom{-}
        {\left.\frac{\partial u{(x_1,x_2)}}{\partial x_{2}}\right|_{x_2=d}}
        + 2i\beta e^{i\hat{\alpha} x_{1} - i\beta d} &\text{ on } \tilde{\Gamma}_{d}, \\[0.2cm]
        {\left[T^{-} u {(\cdot,b)}\right](x_1)} =
        -{\left.\frac{\partial u(x_{1}, x_{2})}{\partial x_{2}}\right|_{x_2=b}} &\text{ on } \tilde{\Gamma}_{b}.
	\end{cases}
\end{equation*}

Multiply both sides of the Helmholtz equation in $C$ by the complex conjugate of $v \in
H_{\hat{\alpha}}^{1}(C)$ and then integrate over $C$. Due to the Green's formula and
the DtN maps, we obtain the variational formulation: Find $u \in H_{\hat{\alpha}}^{1}(C)$
such that
\begin{equation} \label{eq:varEqu}
	a(u, v) = F(v) \text{ for any } v \in H_{\hat{\alpha}}^{1}(C),
\end{equation}
where
\begin{align}
	a(u, v) &\coloneqq \int_{C} \nabla u(x_{1},x_{2}) \cdot \nabla \bar{v}(x_{1},x_{2}) - k^{2}\tilde{q}(x_{1},x_{2}) u(x_{1},x_{2}) \bar{v}(x_{1},x_{2}) \dif x_{1} \dif x_{2}  \label{eq:sesLinForm} \\ 
	& \quad - \int_{\tilde{\Gamma}_{d}} Tu \bar{v} \dif s - \int_{\tilde{\Gamma}_{b}} T^{-}u \bar{v} \dif s, \nonumber \\
	F(v)    &\coloneqq \int_{\tilde{\Gamma}_{d}} -2i\beta e^{i\hat{\alpha} x_{1} - i\beta d} \bar{v} \dif s. \label{eq:linForm}
\end{align}

\begin{theorem} \label{thm:strongelliptic}
	The sesquilinear form $a(\cdot, \cdot):H_{\hat{\alpha}}^{1}(C) \times H_{\hat{\alpha}}^{1}(C)\rightarrow \mathbb{C} $ defined in (\ref{eq:sesLinForm}) is
	strongly elliptic over $H_{\hat{\alpha}}^{1}(C)$.
\end{theorem}

\begin{Proof}
	By the mapping properties of $T^{-}$, we can choose an $N \in \mathbb{N}$ such 
	that, for any $\vert n\vert \geqslant N$ {(cf.~Lemma \ref{lem:unibdd})},
	\begin{equation} \label{eq:negaDtN}
		-\frac{{u_{n}^{-}}'(b)}{u_{n}^{-}(b)} = -\mu_{n, 1} - \frac{v_{n,1}'(b)}{v_{n,1}(b)} < -\mu_{n, 1} + O(M) < 0.
	\end{equation}
	For any quasi-periodic function $f \in H_{\hat{\alpha}}^{1/2}(\mathbb{R})$, defining
	\begin{equation} \label{eq:T2minus}
		T_{1}^{-}f = \sum_{\vert n\vert > N}-\frac{{u_{n}^{-}}'(b)}{u_{n}^{-}(b)}f_{n}e^{i\alpha_{n} x_{1}} \ \text{ and }\
		T_{2}^{-}f = \sum_{\vert n\vert \leqslant N}-\frac{{u_{n}^{-}}'(b)}{u_{n}^{-}(b)}f_{n} e^{i\alpha_{n} x_{1}},
	\end{equation}
    we have $T^{-}f = T_{1}^{-}f + T_{2}^{-}f$. Then we can decompose $a(u, v)$ 
	into two parts such that
	\begin{equation*}
		a(u, v) = B(u, v) - Q(u, v),
	\end{equation*}
	where
	\begin{align*}
		B(u, v) &= \int_{C} \nabla u \cdot \nabla \bar{v} + u\bar{v} \dif x - \int_{\tilde{\Gamma}_{d}} Tu \bar{v} \dif s - \int_{\tilde{\Gamma}_{b}} T_{1}^{-}u \bar{v} \dif s, \\
		Q(u, v) &= \int_{C} \left(k^{2}\tilde{q}(x_{1},x_{2}) + 1\right) u(x_{1},x_{2}) \bar{v}(x_{1},x_{2}) \dif x_{1} \dif x_{2} + \int_{\tilde{\Gamma}_{b}} T_{2}^{-}u \bar{v} \dif s.
	\end{align*}
	
	Taking $v = u$ in $B(u, v)$ gives
	\begin{align*}
		B(u, u) &= \int_{C} \vert\nabla u\vert^{2} + \vert u\vert^{2} \dif x - \int_{\tilde{\Gamma}_{d}} Tu \bar{u} \dif s + {\tt p} \sum_{\vert n\vert > N} \frac{u_{n}^{-'}(b)}{u_{n}^{-}(b)} \vert u_{n}\vert^{2}.
	\end{align*}
	It follows from the definition of $T$ that, for the Fourier expansion 
	$u(x_{1}, d) = \sum_{n} u_{n} e^{i \alpha_{n} x_{1}}$,
    \begin{align*}
        -\Re\left(\int_{\tilde{\Gamma}_{d}} Tu \bar{u} \dif s\right) &= -\Re\left(\int_{\tilde{\Gamma}_{d}}  \left(\sum_{n \in \mathbb{N}} i\beta_{n}u_{n}e^{i\alpha_{n}x_{1}}\right) \left(\sum_{n \in \mathbb{N}} \bar{u}_{n}e^{-i\alpha_{n}x_{1}}\right) \dif s\right) \\
                          &= -\Re\left({\tt p} \sum_{n \in \mathbb{N}} i\beta_{n} \vert u_{n}\vert^{2}\right) \\
                          &= {\tt p} \sum_{\vert \alpha_{n}\vert > k} \vert \beta_{n}\vert \vert u_{n}\vert^{2} \\
                          &\geqslant 0.
    \end{align*}
    By the inequality (\ref{eq:negaDtN}), we have that $\Re B(u, u) \geqslant \Vert u\Vert_{H_{\hat{\alpha}}^{1}(C)}^{2}$. Hence,
	$B(u, v)$ is coercive.

	Next, we claim that $\Re Q(u, v)$ is a compact form. To see this, we note that
	the Riesz representation asserts the existence of two bounded linear operators
	$K_{1} \colon H_{\hat{\alpha}}^{1}(C) \to H_{\hat{\alpha}}^{1}(C)$ and $K_{2}
	\colon H_{\hat{\alpha}}^{1}(C) \to H_{\hat{\alpha}}^{1}(C)$ such that
	\begin{align*}
		(K_{1}u, v)_{H_{\hat{\alpha}}^{1}(C)} &= \Re\int_{C} (k^{2}\tilde{q}(x_{1},x_{2}) + 1) u(x_{1},x_{2}) \bar{v}(x_{1},x_{2}) \dif x_{1} \dif x_{2},\\
		(K_{2}u, v)_{H_{\hat{\alpha}}^{1}(C)} &= \Re\int_{\tilde{\Gamma}_{b}} T_{2}^{-}u \bar{v} \dif s.
	\end{align*}
	Let $v = K_{1}u$, we have that
	\begin{equation} \label{eq:compest}
		\Vert K_{1}u\Vert_{H_{\hat{\alpha}}^{1}(C)}^{2} \leqslant C \Vert u\Vert_{L^{2}(C)} \Vert K_{1}u\Vert_{L^{2}(C)} \leqslant C \Vert u\Vert_{L^{2}(C)} \Vert K_{1}u\Vert_{H_{\hat{\alpha}}^{1}(C)}.
	\end{equation}
	Suppose that $\{u_{n}\}_{n = 1}^{+\infty}$ is a bounded sequence in $H_{\hat{\alpha}}^{1}(C)$.
	Then by the Rellich-Kondrachov theorem, there exists a subsequnce $\{u_{n_{i}}\}_{i
	= 1}^{+\infty}$ such that $\{u_{n_{i}}\}$ is (strongly) convergent in $L^{2}(C)$.
    By the estimate (\ref{eq:compest}), $\{K_{1}u_{n_{i}}\}$ is a Cauchy sequence in
    $H_{\hat{\alpha}}^{1}(C)$ and hence $\{K_{1}u_{n_{i}}\}$ is convergent in $H_{\hat{\alpha}}^{1}(C)$.
    This implies that $K_{1}$ is compact. Since there are only finitely many terms 
	in the summation of $T_{2}^{-}u$ (cf.~(\ref{eq:T2minus}))), it follows that $K_{2}$ 
	is compact.
 	
 	Hence, $\Re a(u, v)$ is the sum of a coercive form and a compact form, implying 
	the strong ellipticity of the sesquilinear form $a(u, v)$ over $H_{\hat{\alpha}}^{1}(C)$.
\end{Proof}

\subsection{Well-posedness results}\label{well-posedness}

In this section we shall prove that our scattering problem admits a unique solution 
for all wavenumbers $k>0$ except for a discrete set. The exceptional wavenumbers in 
this discrete set can accumulate only at infinity and are bounded below by a positive 
constant.  

First, we study the unique solvability at small wavenumbers. Lemma \ref{lem:estimateEtan} 
shows the asymptotic behavior of $\eta_{n}$ as $n$ tends to infinity at a fixed 
wavenumber $k$. Motivated by the arguments in the proof of Lemma \ref{lem:estimateEtan}, 
we derive the behavior of $\eta_{n}$ for fixed $n$ and $k$ tending to $0$. Some 
technical proofs will be presented in the Appendix. 

\begin{lemma} \label{lem:allEtan}
Suppose that $k>0$ is sufficiently small. Then:
\begin{description}
	\item[(i)] $\eta_{n} > 2$ for all $n \neq 0$.
	\item[(ii)] For $n = 0$, we have $\vert \eta_{0}\vert < 2$ if $q(x_{2}) > \sin^{2}\theta$ in $x_2<b$.
\end{description}
\end{lemma}

The proof of the first assertion is given by Lemma \ref{lem:etan_k} in 
Appendix \ref{appendix:asymptotic} and that of the second one in Appendix 
\ref{appendix:discriminant}. By Lemma \ref{lem:allEtan}, we can compute 
the coefficients in the downward DtN map $T^{-}$ explicitly. Recall from \eqref{eq:drc} 
and (\ref{eq:newDtNminus}) that $u_n^-(x_2) e^{i\alpha_{n} x_{1}}$ denotes the $n$-th 
downward wave mode. Using Definition \ref{def:modes} and Lemma \ref{lem:solvereal} (a) 
and (c), we obtain for small wavenumbers that
\begin{equation}
  -\frac{{u_{n}^{-}}'(b)}{u_{n}^{-}(b)} =
  \begin{cases}
      -\mu_{n, 1} - \frac{v_{n,1}'(b)}{v_{n,1}(b)} &\text{ for } n \neq 0, \\
      i\theta_{0} - \frac{v_{0,1}'(b)}{v_{0,1}(b)} &\text{ for } n = 0.
  \end{cases}
\end{equation}
where $\mu_{n,1}>0$ if $n\neq 0$ and $\theta_0\in (0, 1/2)$. Note that $u_n^-$, 
$v_{n,1}$, $\mu_{n, 1}$ and $\theta_0$ all depend on the wavenumber. Their 
asymptotic properties as $k \rightarrow +0$ are shown below.

\begin{lemma} \label{lem:asymptotic} 
Suppose that $k > 0$ is sufficiently small. 
\begin{description} 
    \item[(i)] For $n \neq 0$, the term $\frac{v_{n,1}'(b)}{v_{n,1}(b)}$ is
               real-valued and 
			   \begin{equation*}
			       \left\vert \frac{v_{n,1}'(b)}{v_{n,1}(b)}\right\vert=O(k^2)\quad\mbox{as}\quad k\rightarrow +0.
			   \end{equation*} 
    \item[(ii)] For $n = 0$ and $q(x_{2}) > \sin^{2}\theta$ for $x_{2} < b$. 
	            Then it holds that
                \begin{equation*}
                	\theta_0=O(k),\quad \left\vert \frac{v_{0,1}'(b)}{v_{0,1}(b)}\right\vert=o(1)\quad\mbox{as}\quad k\rightarrow +0.
		        \end{equation*}
    \item[(iii)] For $n \neq 0$, we have $\mu_{n, 1} > \vert n\vert - \frac{1}{2\pi} \ln 2 > 0$ as $k \to +0$. 
\end{description}
\end{lemma}

\begin{Proof}
    For the proof of the asymptotic results in (i), (ii) and (iii) we refer to 
	Lemmas \ref{lem:etan_k}, \ref{lem:unibddk} and \ref{lem:C3} in the  Appendix
    \ref{appendix:asymptotic} (Note that, for similarity reasons, we prove the 
	results for the upper radiation condition instead of the analogous results 
	for the lower radiation condition). It remains to prove that $\frac{v_{n,1}'(b)}{v_{n,1}(b)}$
    is real-valued for $n \neq 0$ if $k > 0$ is sufficiently small. We retain 
	the notation used in Section \ref{sec:radiation}. Recall from  Lemma \ref{lem:solvereal}
    that $v_{n,1}$ is a linear combination of the fundamental solutions of (\ref{eq:hill}), 
	that is, $$v_{n, 1}(x_{2}) =
    e^{-\mu_{n, 1}x_{2}}(c_{1}w_{n, 1}(x_{2}) + c_{2}w_{n, 2}(x_{2})),$$ 
	where $(c_{1}, c_{2})^{\mathrm{T}}$ is an eigenvector associated to the eigenvalue
    $\lambda_{n, 1}$ of the matrix $W_{n}$. Since $q(x_{2})$ is real-valued, the
    fundamental solutions $w_{n, 1}(x_{2})$ and $w_{n, 2}(x_{2})$ are real-valued.
    Due to $\eta_{n} > 2$ (cf.~Lemma \ref{lem:allEtan}), $\mu_{n, 1}$ and $(c_{1},
    c_{2})^{\mathrm{T}}$ are also real-valued. Hence $v_{n, 1}(x_{2})$ is real-valued,
    which completes the proof.
\end{Proof}

\begin{theorem}[Solvability at small wavenumbers] \label{thm:uniqueness}
	Assume that, for any $b - 2\pi < x_{2} \leqslant b$, we have $q(x_{2}) > \sin^{2}\theta$. 
	Then the variational formulation \eqref{eq:varEqu} admits a unique solution 
	$u \in H_{\hat{\alpha}}^{1}(C)$ for all $k\in (0, k_{0})$ and all directions 
	of incidence.
\end{theorem}

\begin{Proof}
    By Theorem \ref{thm:strongelliptic} and the Fredholm alternative,  it is sufficient 
	to prove that the solution of the variational formulation \eqref{eq:varEqu} is unique 
	if $k > 0$ is sufficiently small. Taking $v = u$ in the variational formulation \eqref{eq:varEqu} 
	with $u^{\mathrm{in}}=0$, we see
    \begin{equation} \label{eq:var0}
    	0 = a(u, u) = \int_{C} \vert \nabla u(x)\vert^{2} - k^{2}\tilde{q}(x) \vert u(x)\vert^{2} \dif x - \int_{\tilde{\Gamma}_{d}} Tu \bar{u} \dif s - \int_{\tilde{\Gamma}_{b}} T^{-}u \bar{u} \dif s.
    \end{equation}
    For $k \in (0, k_0)$ with $k_0>0$ sufficiently small, the upward and downward DtN maps $T$ and $T^-$ can be rephrased as
    \begin{align*}
        &\int_{\tilde{\Gamma}_{d}} Tu \bar{u} \dif s = {\tt p} i k\cos\theta \vert u_{0}(d)\vert^{2} - {\tt p} \sum_{n \neq 0} \vert \beta_{n}\vert \vert u_{n}(d)\vert^{2}, \\
        &\int_{\tilde{\Gamma}_{b}} T^{-}u \bar{u} \dif s = {\tt p} \left(i\theta_{0} - \frac{v_{0,1}'(b)}{v_{0,1}(b)}\right)\vert u_{0}(b)\vert^{2} - {\tt p} \sum_{n \neq 0} \left(\mu_{n, 1} + \frac{v_{n, 1}'(b)}{v_{n, 1}'(b)}\right) \vert u_{n}(b)\vert^{2},
    \end{align*}
    where $u_{n}(d)$ and $u_{n}(b)$ are the Fourier coefficients of $e^{-i\alpha
    x_{1}}u(x_{1}, x_{2})$ on $\tilde{\Gamma}_{d}$ and $\tilde{\Gamma}_{d}$, 
	respectively. The real and imaginary parts can be written as
    \begin{align*}
       &\Re \int_{\tilde{\Gamma}_{d}} Tu \bar{u} \dif s = - {\tt p} \sum_{n \neq 0} \vert \beta_{n}\vert \vert u_{n}(d)\vert^{2}<0, \\
       &\Im \int_{\tilde{\Gamma}_{d}} Tu \bar{u} \dif s = {\tt p} k \cos\theta \vert u_{0}(d)\vert^{2}>0, \\	
       &\Re \int_{\tilde{\Gamma}_{b}} T^{-}u \bar{u} \dif s = - {\tt p} \frac{v_{0,1}'(b)}{v_{0,1}(b)}\vert u_{0}(b)\vert^{2} - {\tt p} \sum_{n \neq 0} \left(\mu_{n, 1} + \frac{v_{n, 1}'(b)}{v_{n, 1}'(b)}\right) \vert u_{n}(b)\vert^{2},\\
       &\Im \int_{\tilde{\Gamma}_{b}} T^{-}u \bar{u} \dif s = {\tt p} \theta_{0} \vert u_{0}(b)\vert^{2}>0.
    \end{align*}

    Taking the real part of (\ref{eq:var0}) and using the asymptotic behavior of 
	$\mu_{n, 1}$ and Lemma \ref{lem:asymptotic} (ii), we arrive at
    \begin{equation} \label{eq:realVar}
    \begin{aligned}
        0 & = \Re \,[a(u, u)] \\ 
		  & = \Vert \nabla u\Vert_{L^{2}(C)}^{2} - k^{2} \int_{C} \Re \tilde{q}(x)\vert u(x)\vert^{2} \dif x  - \Re \int_{\tilde{\Gamma}_{d}} Tu \bar{u} \dif s - \Re \int_{\tilde{\Gamma}_{b}} T^{-}u \bar{u} \dif s \\
          & \geqslant \Vert \nabla u\Vert_{L^{2}(C)}^{2} - C_0k^{2} \Vert u\Vert_{L^{2}(C)}^{2} + \frac{{\tt p}}{C_0} \sum_{n \neq 0} \vert n\vert \vert u_{n}(b)\vert^{2} + \frac{{\tt p}}{C_0} \vert u_{0}(b)\vert^{2} - \frac{2\pi}{C_0} \vert u_{0}(b)\vert^{2}+o(1) \vert u_{0}(b)\vert^{2} \\
          & \geqslant C_1 \Vert u\Vert_{H^{1}(C)}^{2} - \frac{{\tt p}}{C_{0}} \vert u_{0}(b)\vert^{2},
    \end{aligned}
    \end{equation}
    where the constants $C_0, C_1 > 0$ do not depend on the direction of incidence. 
    Taking the imaginary part of (\ref{eq:var0}) and letting $k$ tend to zero give that
    \begin{equation} \label{eq:imagVar}
    \begin{aligned}
        0 & = \Im \,[a(u, u)] \\ 
		  & = -k^{2} \int_{C} \Im \tilde{q}(x) \vert u(x)\vert^{2} \dif x - \Im \int_{\tilde{\Gamma}_{d}} Tu \bar{u} \dif s - \Im \int_{\tilde{\Gamma}_{b}} T^{-}u \bar{u} \dif s \\
          & \geqslant -{\tt p} k\cos\theta \vert u_{0}(d)\vert^{2} - {\tt p} \theta_{0}\vert u_{0}(b)\vert^{2}.
    \end{aligned}
    \end{equation}
    Since $\cos\theta$ and $\theta_{0}$ are positive, we have that $\vert u_{0}(d)\vert
    = \vert u_{0}(b)\vert = 0$. Combining this and \eqref{eq:realVar} implies 
	that $\Vert u\Vert_{H^{1}(C)} = 0$, which completes the uniqueness proof 
	for small wavenumbers. 
\end{Proof}

\begin{theorem}\label{mainth}
	Fix the angle of incidence $\theta \in (-\pi/2, \pi/2)$. Then the variational
    problem \eqref{eq:varEqu} has a unique solution for all wavenumbers $k > 0$ not
    contained in a  discrete subset of the positive reals. The only possible 
	accumulating point of this discrete set is infinity.
\end{theorem}

\begin{Proof}
	Without loss of generality we may assume that, for any $b-2\pi< x_2 \leqslant b$, 
	we have $q(x_{2}) > \sin^{2}\theta$. Otherwise we could replace $q$ by $q_{C} \coloneqq C^2q$ 
	and $k$ by $k_C \coloneqq k/C$ with a big constant $C>0$. Then our assumption is fulfilled
    with $q_C$ and $k_C$, and the resulting theorem for $q_C$ and $k_C$ implies 
	the theorem for $q$ and $k$.

    We shall carry out the proof by applying the analytic Fredholm theory 
	\cite[Theorem 8.26]{CoKr} to our wavenumber-dependent variational formulation. 
	For this purpose we need to transfer the dependence of the $\alpha$-quasiperodic 
	Sobolev space on $k$ to an equivalent variational formulation, for which the 
	underlying function space is independent of the wavenumber. Setting $\hat{a}(u, v) 
	\coloneqq a(\hat{u},\hat{v})$ for $u, v \in H^1_{per}(C)$ with $\hat{u} \coloneqq 
	e^{-i \hat{\alpha} x_{1}}u$ as well as $\hat{v} \coloneqq e^{-i \hat{\alpha} x_{1}}v$, 
	we obtain a new sesquilinear form $\hat{a}(\cdot, \cdot)$ defined on $H_{per}^{1}(C)$:
    \begin{equation} \label{eq:pVarE}
		\begin{aligned}
			\hat{a}(u, v) = &\int_{C} \nabla u(x) \cdot \nabla \bar{v}(x) - i \hat{\alpha} \frac{\partial u(x)}{\partial x_{1}}\bar{v}(x) + i \hat{\alpha} u(x) \frac{\partial \bar{v}(x)}{\partial x_{1}} - (k^{2}\tilde{q}(x) - \hat{\alpha}^{2})u(x)\bar{v}(x)\dif x \\
			                & - \int_{\tilde{\Gamma}_{d}} \hat{T}u \bar{v} \dif s - \int_{\tilde{\Gamma}_{b}} \hat{T}^{-}u \bar{v} \dif s.
		\end{aligned}
    \end{equation}
    Here the operators $\hat{T}$ and $\hat{T}^{-}$ are periodic DtN maps defined over 
	the space $H_{per}^{1/2}(\tilde\Gamma_d)$ and $H_{per}^{1/2}(\tilde\Gamma_b)$, 
	respectively. Using the decomposition of $T^{-}$ in (\ref{eq:T2minus}), we 
	can rewrite $\hat{a}(u, v)$ as $\hat{a}_{1}(u, v) - \hat{a}_{2}(u, v)$, where
    \begin{align*}
        \hat{a}_{1}(u, v) &\coloneqq \int_{C} \nabla u \cdot \nabla \bar{v} + u \bar{v} \dif x - \int_{\tilde{\Gamma}_{d}} \hat{T}u \bar{v} \dif s - \int_{\tilde{\Gamma}_{b}} \hat{T}_{1}^{-}u \bar{v} \dif s, \\
        \hat{a}_{2}(u, v) &\coloneqq \int_{C} i \hat{\alpha} \frac{\partial u(x)}{\partial x_{1}}\bar{v}(x) - i \hat{\alpha} u(x) \frac{\partial \bar{v}(x)}{\partial x_{1}} + (k^{2}\tilde{q}(x) - \hat{\alpha}^{2} + 1)u(x)\bar{v}(x)\dif x + \int_{\tilde{\Gamma}_{b}} \hat{T}_{2}^{-}u \bar{v} \dif s.
    \end{align*}
    By the Riesz representation theorem, the sesquilinear forms $\hat{a}_{1}(\cdot, \cdot)$ 
	and $\hat{a}_{2}(\cdot, \cdot)$ give rise to the wavenumber-dependent operators 
	$A(k)$ and $K(k)$, respectively, where $A(k)$ is a bounded invertible operator 
	on $H_{per}^{1}(C)$ and $K(k): H_{per}^{1}(C) \to H_{per}^{1}(C)$ is compact 
	(cf.~the proof of Theorem \ref{thm:strongelliptic} and \cite[Section 4]{Be}). 
	It follows from \cite[Lemma 6]{Ki} and the analytic dependence of the wave modes 
	on $k$ that both, $A(k)$ and $K(k)$, are analytic operators with respect to the 
	wavenumber $k$. Combining the uniqueness result (cf.~Theorem \ref{thm:uniqueness}) 
	and the argument in \cite[Theorem 2.1]{Ba2}, we can show that $(A(k) - K(k))^{-1}$
    exists for all $k$ except those in a discrete set. By the analytic Fredholm theory, 
    the only possible accumulating point of this discrete set is at infinity. 
    This concludes the proof.
\end{Proof}


\begin{appendices}
\section{Floquet theory} \label{appendix:floquet}
This appendix is mainly based on \cite[Chapter 1]{Ea2} and \cite[Chapter 9]{Ro}. 
Consider the general second-order ODE
\begin{equation} \label{eq:generalODE}
	a_{0}(x) y''(x) + a_{1}(x) y'(x) + a_{2}(x) y(x) = 0,
\end{equation}
where the coefficients $a_{0}(x)$, $a_{1}(x)$ and $a_{2}(x)$
are complex-valued, piecewise continuous and periodic with the period $p$.

\begin{theorem} \label{thm:property}
	There are a non-zero constant $\rho$ and a non-trival solution $\psi(x)$ of
	(\ref{eq:generalODE}) such that
    \begin{equation} \label{eq:property}
		\psi(x + p) = \rho \psi(x).
	\end{equation}
\end{theorem}

\begin{Proof}
	Let $\phi_{1}(x)$ and $\phi_{2}(x)$ be the linearly independent solutions
	of (\ref{eq:generalODE}) which satisfy the initial conditions
	\begin{equation} \label{eq:initalcondition}
		\begin{aligned}
			\phi_{1}(0) = 1,\ \phi_{1}'(0) = 0; \\
			\phi_{2}(0) = 0,\ \phi_{2}'(0) = 1.
		\end{aligned}
	\end{equation}
	Since $\phi_{1}(x + p)$ and $\phi_{2}(x + p)$ are also linearly independent
	solution of (\ref{eq:generalODE}), there are constants $A_{ij}$($1 \leqslant i,
	j \leqslant 2$) such that
	\begin{equation} \label{eq:afterperiod}
		\begin{aligned}
			\phi_{1}(x + p) = A_{11}\phi_{1}(x) + A_{12}\phi_{2}(x), \\
			\phi_{2}(x + p) = A_{21}\phi_{1}(x) + A_{22}\phi_{2}(x),
		\end{aligned}
	\end{equation}
	where the matrix $A = (A_{ij})$ is non-singular. Every solution $\psi(x)$ of
	(\ref{eq:generalODE}) has the form
	\begin{equation} \label{eq:psi_func}
		\psi(x) = c_{1}\phi_{1}(x) + c_{2}\phi_{2}(x),
	\end{equation}
	where $c_{1}$ and $c_{2}$ are constant. By (\ref{eq:afterperiod}), we have
	that (\ref{eq:property}) holds if
	\begin{equation*} 
		\begin{aligned}
			(A_{11} - \rho) c_{1} + A_{21} c_{2} = 0, \\
			A_{12} c_{1} + (A_{22} - \rho) c_{2} = 0.
		\end{aligned}
	\end{equation*}
    The function $\psi$ of (\ref{eq:psi_func}) is non-trival if $c_{1}$ and $c_{2}$ are not both zero. Hence, a non-trival $\psi$ exists if and only if
    \begin{equation*}
    	\Det\begin{pmatrix}
    		A_{11} - \rho & A_{21} \\
    		A_{12}        & A_{22} - \rho
    	\end{pmatrix} = 0,
    \end{equation*}
	i.e.,
	\begin{equation} \label{eq:quadratic}
		\rho^{2} - (A_{11} + A_{22}) \rho + \Det A = 0.
	\end{equation}
	This is a quadratic equation for $\rho$. Since $A$ is non-singular, (\ref{eq:quadratic})
	has at least one non-trival solution. This completes the proof.
\end{Proof}

An alternative form of (\ref{eq:quadratic}) can be obtained as follows. It follows
from (\ref{eq:initalcondition}) and (\ref{eq:afterperiod}) that
\begin{equation*}
	A_{11} = \phi_{1}(p),\ A_{12} = \phi_{1}'(p),\ A_{21} = \phi_{2}(p),\ A_{22} = \phi_{2}'(p).
\end{equation*}
Hence using Liouville's formula for the Wronskian
\begin{equation*}
    W(\phi_{1}, \phi_{2})(x) \coloneqq \Det \begin{pmatrix}
    		 \phi_{1}(x)&\phi_{2}(x)\\
    		\phi_{1}'(x)&\phi_{2}'(x)
    	\end{pmatrix}
\end{equation*}
\cite[Lemma 3.11]{Te} and the fact that $W(\phi_{1}, \phi_{2})(0) = 1$, we get
\begin{equation*}
	\Det A = W(\phi_{1}, \phi_{2})(p) = e^{\left(-\int_{0}^{p} \frac{a_{1}(x)}{a_{0}(x)} \dif x\right)}.
\end{equation*}
Thus (\ref{eq:quadratic}) can be rewritten as
\begin{equation} \label{eq:quadratic2}
    \begin{aligned}
        &\rho^{2} - \eta \rho + e^{\left(-\int_{0}^{p} \frac{a_{1}(x)}{a_{0}(x)} \dif x\right)} = 0, \\
        &\eta = \phi_{1}(p) + \phi_{2}'(p).
    \end{aligned}
\end{equation}
In particular, all of the zeros $\varrho$ of (\ref{eq:quadratic2}) or equivalently of (\ref{eq:quadratic}) are different from zero.

Theorem \ref{thm:property} leads to the existence of two linearly independent
solutions of (\ref{eq:generalODE}) having the special form given in Theorem
\ref{thm:floquet} below.

\begin{theorem} \label{thm:floquet}
	There are linearly independent solutions $\psi_{1}(x)$ and $\psi_{2}(x)$ of
	(\ref{eq:generalODE}) such that either
	\begin{enumerate}[(1)]
		\item (\ref{eq:quadratic2}) has two distinct roots or only one root $\rho$ 
		with $\Rank(A^{\mathrm{T}} - \rho I) = 0$.
        \begin{equation*}
			\psi_{1}(x) = e^{m_{1}x} v_{1}(x),\ \psi_{2}(x) = e^{m_{2}x} v_{2}(x),
		\end{equation*}
		where $m_{1}, m_{2}$ are constant, not necessarily distinct, and $v_{1},
		v_{2}$ are periodic functions with period $p$, or
		\item (\ref{eq:quadratic2}) has only one root $\rho$ with $\Rank(A^{\mathrm{T}} - \rho I) = 1$.
        \begin{equation*}
			\psi_{1}(x) = e^{mx} v_{1}(x),\ \psi_{2}(x) = e^{mx}(x v_{1}(x) + v_{2}(x)),
		\end{equation*}
		where $m$ is constant, and $v_{1}, v_{2}$ are periodic with period $p$.
	\end{enumerate}
\end{theorem}

\begin{Proof}
	Consider two cases as follows.
	\begin{enumerate}[(1)]
		\item Suppose that (\ref{eq:quadratic}) has distinct solutions $\rho_{1}$
		      and $\rho_{2}$. Then by the proof of Theorem \ref{thm:property}, 
			  there are non-trival solutions $\psi_{1}$ and $\psi_{2}$ of (\ref{eq:generalODE}) 
			  such that
		      \begin{equation} \label{eq:property2}
		      	\psi_{k}(x + p) = \rho_{k} \psi_{k}(x),\ k = 1, 2.
		      \end{equation}
	          It is easy to show that $\psi_{1}$ and $\psi_{2}$ are linearly
	          independent. Since $\rho_{1}$ and $\rho_{2}$ are non-zero, we can
	          define $m_{1}$ and $m_{2}$ so that
	          \begin{equation} \label{eq:exponents}
	          	e^{p m_{k}} = \rho_{k},\ k = 1, 2.
	          \end{equation}
              We can construct a new function $v_{k}(x)$,
              \begin{equation}
              	v_{k}(x) = e^{-m_{k}x} \psi_{k}(x).
              \end{equation}
          	  By (\ref{eq:property2}) and (\ref{eq:exponents}),
          	  \begin{equation*}
          	  	v_{k}(x + p) = e^{-m_{k}(x + p)} \rho_{k} \psi_{k}(x) = v_{k}(x).
          	  \end{equation*}
              Thus,
              \begin{equation*}
              	\psi_{k}(x) = e^{m_{k}x} v_{k}(x),
              \end{equation*}
              where $v_{k}(x)$, $k = 1, 2$ are periodic functions.
		\item Suppose that (\ref{eq:quadratic}) has a repeated solution $\rho$.
		      We define $m$ so that
		      \begin{equation*}
		      	e^{pm} = \rho.
		      \end{equation*}
	          By Theorem \ref{thm:property}, there is a non-trival solution $\Psi_{1}$
	          of (\ref{eq:generalODE}) such that
	          \begin{equation} \label{eq:property3}
	          	\Psi_{1}(x + p) = \rho \Psi_{1}(x).
	          \end{equation}
              Let $\Psi_{2}$ be any solution of (\ref{eq:generalODE}) which is
              linearly independent of $\Psi_{1}$. Since $\Psi_{2}$ also satisfies
              (\ref{eq:generalODE}), there are constants $d_{1}$ and $d_{2}$ such that
              \begin{equation} \label{eq:afterperiod2}
              	\Psi_{2}(x + p) = d_{1}\Psi_{1}(x) + d_{2}\Psi_{2}(x)
              \end{equation}
              Next we calculate $d_{2}$. With the help of (\ref{eq:property3}) and (\ref{eq:afterperiod2}), we obtain
              \begin{equation*}
              	W(\Psi_{1}, \Psi_{2})(x + p) = \rho d_{2} W(\Psi_{1}, \Psi_{2})(x).
              \end{equation*}
              Hence by Liouville's formula for the Wronskian,
              \begin{equation*}
              	\rho d_{2} = e^{-\int_{x}^{x+p} \frac{a_{1}(t)}{a_{0}(t)} \dif t} = e^{-\int_{0}^{p} \frac{a_{1}(t)}{a_{0}(t)} \dif t},
              \end{equation*}
              since the integrand has period $p$. Futhermore the right-hand side
              in the above equation is equal to $\rho^{2}$ since $\rho$ is the repeated
              solution of (\ref{eq:quadratic2}). Hence $d_{2} = \rho$. From (\ref{eq:afterperiod2}) we now have
              \begin{equation} \label{eq:afterperiod3}
              	\Psi_{2}(x + p) = d_{1}\Psi_{1}(x) + \rho \Psi_{2}(x).
              \end{equation}

              Now, there are two possibilities to consider:

              Possibility 1: $d_{1} = 0$

              In this case (\ref{eq:afterperiod3}) implies
              \begin{equation*}
              	\Psi_{2}(x + p) = \rho \Psi_{2}(x).
              \end{equation*}
          	  This together with (\ref{eq:property3}) indicates that we have the
          	  same situation as in (\ref{eq:property2}) but with $\rho_{1} = \rho_{2}
          	  = \rho$. In other words, $(A^{\mathrm{T}} - \rho I) \mathbf{x} = 0$ has two linearly independent solutions, i.e., $\Rank(A^{\mathrm{T}} - \rho I) = 0$. Consequently, this case is covered by part (1) of the theorem.

              Possibility 2: $d_{1} \neq 0$

              In this case we define
              \begin{align*}
              	v_{1}(x) &\coloneqq e^{-mx} \Psi_{1}(x),\\
              	v_{2}(x) &\coloneqq e^{-mx} \Psi_{2}(x) - \left(\frac{d_{1}}{p \rho}\right) x v_{1}(x).
              \end{align*}
              Then due to (\ref{eq:property3}) and (\ref{eq:afterperiod3}) $v_{1}$
              and $v_{2}$ have period $p$. In this case, $(A^{\mathrm{T}} - \rho I) 
			  \mathbf{x} = 0$ has only one non-trival solution, i.e., $\Rank(A^{\mathrm{T}} - \rho I) = 1$. 
			  Therefore, since
              \begin{align*}
              	\Psi_{1}(x) &= e^{mx}v_{1}(x), \\
              	\Psi_{2}(x) &= e^{mx}\left(\left(\frac{d_{1}}{p \rho}\right) x v_{1}(x) + v_{2}(x)\right),
              \end{align*}
              part (2) is covered with $\psi_{1}(x) \coloneqq \Psi_{1}(x)$ and $\psi_{2}(x)
              \coloneqq (p\rho / d_{1}) \Psi_{2}(x)$.
	\end{enumerate}
\end{Proof}

\section{Proof of Lemma \ref{lem:allEtan} (ii)} \label{appendix:discriminant}

In this section, we mainly follow \cite[Section 2.1 and 2.2]{Ea2}. For simplicity,
we omit the proofs presented in this reference paper. Here we consider the 
second-order equation
\begin{equation} \label{eq:evp}
    (p(x)y'(x))' + (\lambda s(x) - q(x))y(x) = 0,
\end{equation}
where $\lambda$ is a real parameter. The coefficient functions $p(x)$, $s(x)$ 
and $q(x)$ are real-valued and periodic with period $p$. Furthermore, we assume 
the existence of a constant $s > 0$ such that $s(x) \geqslant s$.

To show the dependence of the solutions on $\lambda$, we denote the fundamental 
solutions of (\ref{eq:evp}) by $\phi_{1}(x, \lambda)$ and $\phi_{2}(x, \lambda)$. 
Then corresponding to (\ref{eq:quadratic2}), we define
\begin{equation*}
    \eta(\lambda) \coloneqq \phi_{1}(p, \lambda) + \phi_{2}'(p, \lambda).
\end{equation*}
Whether $\lambda$ is real or complex, $\phi_{1}(x, \lambda)$ and $\phi_{2}(x, \lambda)$
and their derivatives are analytic functions of $\lambda$ for fixed $x$, see 
\cite[Section 1.7]{Ea1}. Hence, $\eta(\lambda)$ is an analytic function of $\lambda$.

Now we introduce two eigenvalue problems associated to (\ref{eq:evp}) on the
interval $[0, p]$, where $\lambda$ is regarded as a parameter. These two problems
will be used in the investigation of $\eta(\lambda)$.

\textbf{The periodic eigenvalue problem}: Consider (\ref{eq:evp}) in $[0, p]$
and periodic boundary conditions:
\begin{equation} \label{eq:pevp}
    \begin{aligned}
        &(p(x)y'(x))' + (\lambda s(x) - q(x))y(x) = 0 \text{ in } [0, p], \\
        &y(0) = y(p),\ y'(0) = y'(p).
    \end{aligned}
\end{equation}
It is a self-adjoint problem and has a countable set of eigenvalues $\lambda$. We 
shall denote the eigenfunctions throughout by $\psi_{n}(x)$ and the eigenvalues by
$\lambda_{n}\ (n = 0, 1, \ldots)$, where
\begin{equation*}
    \lambda_{0} \leqslant \lambda_{1} \leqslant \lambda_{2} \leqslant \cdots, \text{ and } \lambda_{n} \to \infty \text{ as } n \to \infty.
\end{equation*}
The $\psi_{n}(x)$ can be chosen to be real-valued and to form an orthonormal
basis over $[0, p]$ with weight function $s(x)$. Thus
\begin{equation*}
    \int_{0}^{p} \psi_{m}(x)\psi_{n}(x)s(x) \dif x =
    \begin{cases}
        1 &(m = n), \\
        0 &(m \neq n).
    \end{cases}
\end{equation*}

\textbf{The semi-periodic eigenvalue problem}: Consider (\ref{eq:evp}) in $[0, p]$
and semi-periodic boundary conditions:
\begin{equation} \label{eq:spevp}
    \begin{aligned}
        &(p(x)y'(x))' + (\lambda s(x) - q(x))y(x) = 0 \text{ in } [0, p], \\
        &y(0) = -y(p),\ y'(0) = -y'(p).
    \end{aligned}
\end{equation}
It is also a self-adjoint problem and we shall denote the eigenfunctions by
$\xi_{n}(x)$ and the eigenvalues $\lambda$ by $\mu_{n}\ (n = 0, 1, \cdots)$, where
\begin{equation*}
  \mu_{0} \leqslant \mu_{1} \leqslant \mu_{2} \leqslant \cdots, \text{ and } \mu_{n} \to \infty \text{ as } n \to \infty.
\end{equation*}

\begin{figure}[htp]
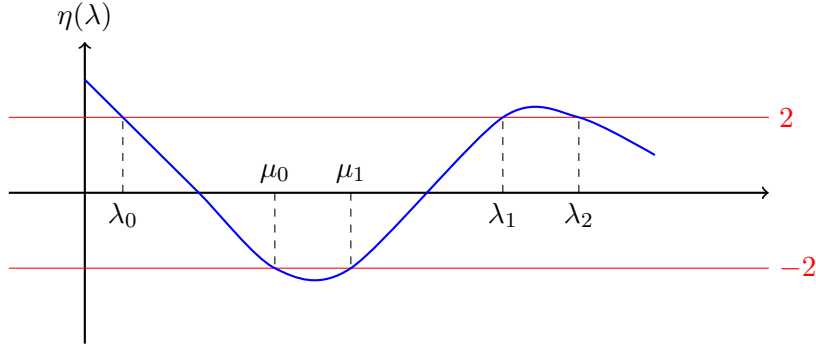

	\centering
	\includestandalone{figs/discriminant}
	\caption{An illustration of dependence of $\eta$ on $\lambda$.}
	\label{fig:discriminant}
\end{figure}

Now we use the existence of the eigenvalues $\lambda_{n}$ and $\mu_{n}$ in the
eigenvalue problems (\ref{eq:pevp}) and (\ref{eq:spevp}) to study the behavior
of $\eta(\lambda)$. We refer to Figure \ref{fig:discriminant} for a rough idea
of the function $\eta(\lambda)$.

\begin{theorem} \label{thm:discriminant}
    \begin{enumerate}[(1)]
      \item The numbers $\lambda_{n}$ and $\mu_{n}$ occur in the order
            \begin{equation*}
                \lambda_{0} < \mu_{0} \leqslant \mu_{1} < \lambda_{1} \leqslant \lambda_{2} < \mu_{2} \leqslant \mu_{3} < \lambda_{3} \leqslant \lambda_{4} < \cdots.
            \end{equation*}
      \item $\eta(\lambda)$ decreases from $2$ to $-2$ in the intervals $[\lambda_{2m}, \mu_{2m}]$.
      \item $\eta(\lambda)$ increases from $-2$ to $2$ in the intervals $[\mu_{2m+1}, \lambda_{2m+1}]$.
      \item $\eta(\lambda) > 2$ in the intervals $(-\infty, \lambda_{0})$ and $(\lambda_{2m+1}, \lambda_{2m+2})$.
      \item $\eta(\lambda) < -2$ in the intervals $(\mu_{2m}, \mu_{2m+1})$.
    \end{enumerate}
\end{theorem}

{\bf{Proof of Lemma \ref{lem:allEtan} (ii). }}

For $n = 0$, the ordinary differential equation (\ref{eq:hill}) takes the form
\begin{equation} \label{eq:hilln0}
    u_{0}''(x_{2}) + k^{2}(q(x_{2}) - \sin^{2}\theta)u_{0}(x_{2}) = 0.
\end{equation}
Then we consider the periodic eigenvalue problem
\begin{equation} \label{eq:pevpn0}
    \begin{aligned}
        &u_{0}''(x_{2}) + \lambda(q(x_{2}) - \sin^{2}\theta)u_{0}(x_{2}) = 0 \text{ in } (0, 2\pi), \\
        &u_{0}(0) = u_{0}(2\pi),\ u_{0}'(0) = u_{0}'(2\pi).
    \end{aligned}
\end{equation}
Integrating by parts and using $q(x_{2}) > \sin^{2}\theta$, one can show that 
the periodic eigenvalue problem (\ref{eq:pevpn0}) has nonnegative eigenvalues. 
Furthermore, $\lambda = 0$ is the first eigenvalue of (\ref{eq:pevpn0}). Applying
Theorem \ref{thm:discriminant} with $\lambda=k^2$, it follows that the constant 
$\eta_{0}$ associated to (\ref{eq:hilln0}) (i.e., $\eta_n$ defined by \eqref{eq:etan} 
with $n=0$) must lie in $(-2, 2)$ if $k > 0$ is sufficiently small.  

\section{Asymptotic behavior of $\frac{v_{n,1}'(d)}{v_{n,1}(d)}$ for $k \!\to\! +0$} \label{appendix:asymptotic}

\begin{lemma} \label{lem:etan_k}
    If the positive wavenumber  $k$ is sufficiently small and if $n \neq 0$, then 
	\begin{equation*}
		\eta_{n} > 2, \quad \mu_{n, 1} > \vert n\vert - \frac{1}{2\pi} \ln 2 > 0\quad\mbox{if}\quad n \neq 0.	
	\end{equation*}
\end{lemma}

\begin{Proof}
    With the same method as in the proof of Lemma \ref{lem:estimateEtan}, we can
    estimate $\eta_{n}$ by constructing subsolutions and supersolutions for ODEs. Consequently,
	\begin{align*}
		\eta_{n} &= w_{n,1}(2\pi) + w_{n,2}'(2\pi) \\
		         &\geqslant \frac{1}{2}e^{2\pi\beta_{n}^{+}} + \frac{1}{2}e^{-2\pi\beta_{n}^{+}} + \frac{1}{2}e^{2\pi\beta_{n}^{+}} + \frac{1}{2}e^{-2\pi\beta_{n}^{+}} \\
		         &> 2.
	\end{align*}
    Furthermore, $\eta_{n} \to e^{2\pi\vert n\vert} + e^{-2\pi\vert n\vert}$ as 
	$k \to +0$. From 
	\begin{equation*}
		\lambda_{n, j}^{2} - \eta_{n} \lambda_{n, j} + 1 = 0, \quad j = 1, 2, 
	\end{equation*}
    we know $\lambda_{n, 1} > \frac{1}{2}(e^{2\pi\vert n\vert}+e^{-2\pi\vert
    n\vert}) > 1$ and $0 < \lambda_{n,2} < 1$ as $k \to +0$. Then it follows that
	\begin{equation*}
		e^{2\pi\mu_{n, 1}} = \lambda_{n, 1} > \frac{1}{2}(e^{2\pi\vert n\vert} + e^{-2\pi\vert n\vert}) > \frac{1}{2}e^{2\pi\vert n\vert}
	\end{equation*}
    as $k \to +0$. Hence, one obtains $\mu_{n, 1} > \vert n\vert - \frac{1}{2\pi} \ln 2
    > 0$ for any $n \neq 0$ as $k \to +0$.
\end{Proof}

In what follows, for simplicity of notation, we drop the subscript $n$ in $\alpha_{n}$, 
$u_n$ and $\mu_{n, j}$, $\lambda_{n,j}$ for simplicity.

\begin{lemma} \label{lem:unibddk}
	For fixed $n \neq 0$, the term $\left\vert \frac{v_{n,1}'(d)}{v_{n,1}(d)}\right\vert=O(k^2)$ as $k \to +0$. 
\end{lemma}

\begin{Proof}
By using arguments similar to the proof of Lemma \ref{lem:unibdd}, we investigate
the asymptotic behavior of $\frac{v_{n,1}'(d)}{v_{n,1}(d)}$ as $k \to +0$ for
fixed $n \neq 0$. The proof consists of three steps.

\textbf{Step 1: Fundamental solutions.} We look for solutions to the differential equation
$u''(t) + (k^{2}q(t) - \alpha^{2}) u(t) = 0$ in the form $u(t) = \psi(t)\varphi(t)$.
Here the $\psi$ are the solutions of $u''(t) - \alpha^{2}u(t) = 0$ with initial conditions
(\ref{eq:initial1}) or (\ref{eq:initial2}). Clearly, we find that $\varphi(t)$
should satisfy the equation
\begin{equation} \label{eq:varphik}
	\varphi(t) = 1 - k^{2} \int_{0}^{t} \Big[q(\sigma)V^{2}(\sigma)(\tilde{H}(t) - \tilde{H}(\sigma))\Big] \varphi(\sigma) \dif\sigma
\end{equation}
where
\begin{equation*}
	\tilde{H}(t) \coloneqq \frac{1}{\alpha}
	\begin{cases}
		\tanh(\alpha t)&\text{ with condition } (\ref{eq:initial1}), \\
		\coth(\alpha t)&\text{ with condition } (\ref{eq:initial2}),
	\end{cases}
\end{equation*}
and
\begin{equation*}
	V(t) \coloneqq
	\begin{cases}
		\cosh(\alpha t)&\text{ with condition } (\ref{eq:initial1}), \\
		\sinh(\alpha t)&\text{ with condition } (\ref{eq:initial2}).
	\end{cases}
\end{equation*}

We can reformulate (\ref{eq:varphik}) as the operator equation
\begin{equation} \label{eq:opequk}
	\varphi(t) + k^{2}K\varphi(t) = 1,\ j = 1, 2,
\end{equation}
where $(K\varphi)(t) = \int_{0}^{t} k(t, \sigma)\varphi(\sigma) \dif\sigma$ and
$k(t, \sigma) = q(\sigma)V^{2}(\sigma)(\tilde{H}(t) - \tilde{H}(\sigma))$. We
note that the kernels $k(t, \sigma)$ are of different form under the different initial
conditions (\ref{eq:initial1}) and (\ref{eq:initial2}). With the estimates of $k(t,\sigma)$, we have that
\begin{equation*}
	\Vert k^{2}K\Vert \leqslant \frac{Ck^{2}}{\alpha} = C\frac{k^{2}}{k\sin\theta + n} = O(k).
\end{equation*}
Since sufficiently small $k$ guarantees that $\Vert k^{2}K\Vert < 1$, by a Neumann
series argument, $(I+k^{2}K)^{-1}$ exists and
\begin{equation} \label{eq:expanPhik}
	\varphi(t) = 1 - k^{2}(K1)(t) + O(k^{2}).
\end{equation}

We also need to estimate $\varphi'(t)$. By (\ref{eq:opequk}) and (\ref{eq:expanPhik})
we have that
\begin{align*}
    \varphi(t) &= 1 - k^{2}(K1)(t) + k^{4}K(K1)(t) + \cdots, \\
    \varphi'(t) &= -[\partial_{t}k^{2}K]\{1 + O(k^{2})\}.
\end{align*}
It remains to estimate $\Vert \partial_{t}k^{2}K\Vert$. From the definition
of $K$,
\begin{equation*}
    \partial_{t}[k^{2}(K\varphi)(t)] = k^{2}\int_{0}^{t} \frac{V^{2}(\sigma)}{V^{2}(t)}q(\sigma)\varphi(\sigma) \dif\sigma.
\end{equation*}
This implies that $\Vert \partial_{t}k^{2}K\Vert \leqslant Ck^{2}$.

\textbf{Step 2: Characteristic exponents.} Using the fundamental solutions in 
Step 1, the definition (\ref{eq:etan}), and $k \rightarrow 0$, we have that
\begin{align*}
	\eta = &\ \frac{1}{2}(e^{2\pi\alpha_n} + e^{-2\pi\alpha_{n}})(1 - O(k)) + \frac{1}{2}(e^{2\pi\alpha_{n}} + e^{-2\pi\alpha_{n}})(1 - O(k)) \\
	       & + \frac{1}{2}(e^{2\pi\alpha_{n}} - e^{-2\pi\alpha_{n}})O(k^{2}),
\end{align*}
 which implies $\lim_{k \to +0} \eta
= e^{2\pi n} + e^{-2\pi n}$ for $\vert n\vert \neq 0$. Then by (\ref{eq:polynomial}),
\begin{align*}
	\lim_{k \to +0}\lambda_{1} &= \lim_{k \to +0}\left[\frac{\eta}{2} + \sqrt{\frac{\eta^{2}}{4} - 1}\; \right]\\
	            &= \frac{1}{2}(e^{2\pi n} + e^{-2\pi n}) + \sqrt{\frac{e^{4\pi n} + 2 + e^{-4\pi n}}{4} - 1} \\
	            &= e^{2\pi n}.
\end{align*}
It follows from the Definition \ref{def:characteristic} that
\begin{align*}
	&\lim_{k \to +0}\mu_{1} = \lim_{k \to +0}\frac{1}{2\pi} \ln \lambda_{1} = \lim_{k \to +0} \alpha, \\
	&\lim_{k \to +0} \mu_{2} = -\lim_{k \to +0} \mu_{1} = -\lim_{k \to +0} \alpha.
\end{align*}
Furthermore, we know for $\vert n\vert \neq 0$ that
\begin{equation} \label{eq:Eq}
	\lim_{k \to +0} [\mu_1 - \alpha] = 0,\quad \lim_{k \to +0} [\mu_2 + \alpha] = 0.
\end{equation}

\textbf{Step 3: Asymptotics of Hill's equation.}
According to the downward modes in case (a) of Definition \ref{def:modes}, we
assume that $u = e^{\mu t}v(t)$, where $v(0) = v(2\pi)$ and $v'(0) =
v'(2\pi)$. Substituting the product presentation of $u$ into $u''(t) + (k^{2}q(t) - \alpha^{2})u(t)
= 0$, we have that
\begin{equation*}
    e^{\mu t}(v''(t)+2\mu v'(t) + (k^{2}q(t)+\mu^{2}-\alpha^{2})v(t)) = 0.
\end{equation*}
Then $(e^{2\mu t}v'(t))' = -e^{2\mu t}(k^{2}q(t)+\mu^{2}-\alpha^{2})v(t)$. Hence
\begin{equation*}
    e^{2\mu t}v'(t) = v'(0) - \int_{0}^{t} e^{2\mu\tau}(k^{2}q(\tau)+\mu^{2}- \alpha^{2})v(\tau) \dif\tau.
\end{equation*}

By straightforward computation, we get the two identities
\begin{equation} \label{eq:derivativeV1k}
    e^{4\pi\mu}e^{2\mu t} v'(t) = e^{4\pi\mu}v'(0) - e^{4\pi\mu}\int_{0}^{t} e^{2\mu\tau}(k^{2}q(\tau)+\mu^{2}-\alpha^{2})v(\tau) \dif\tau
\end{equation}
and
\begin{equation} \label{eq:derivativeV2k}
    e^{2\mu t}v'(t) = e^{4\pi\mu}v'(0) + \int_{t}^{2\pi} e^{2\mu\tau}(k^{2}q(\tau)+\mu^{2}-\alpha^{2})v(\tau) \dif\tau.
\end{equation}
Subtracting (\ref{eq:derivativeV1k}) from (\ref{eq:derivativeV2k}) gives
\begin{align*}
    v'(t) =\  & \frac{1}{1 - e^{4\pi\mu}} \int_{t}^{2\pi} e^{2\mu(\tau - t)}(k^{2}q(\tau) + \mu^{2} - \alpha^{2})v(\tau) \dif\tau \\
           &+\frac{e^{4\pi\mu}}{1 - e^{4\pi\mu}} \int_{0}^{t} e^{2\mu(\tau - t)}(k^{2}q(\tau) + \mu^{2} - \alpha^{2})v(\tau) \dif\tau.
\end{align*}

There are two cases to consider. We first suppose that $v(0) \neq 0$. Without loss
of generality, let $v(0) = 1$. Then
\begin{align*}
    v(t) = &\ 1 + \int_{0}^{t}v'(\sigma) \dif\sigma \\
         = &\ 1 + \frac{1}{1 - e^{4\pi\mu}} \int_{0}^{t}\int_{\sigma}^{2\pi} e^{2\mu(\tau - \sigma)}(k^{2}q(\tau) + \mu^{2} - \alpha^{2})v(\tau) \dif\tau \dif\sigma \\
           &  + \frac{e^{4\pi\mu}}{1 - e^{4\pi\mu}} \int_{0}^{t}\int_{0}^{\sigma} e^{2\mu(\tau - \sigma)}(k^{2}q(\tau) + \mu^{2} - \alpha^{2})v(\tau) \dif\tau \dif\sigma \\
         = &\ 1 + \frac{1}{1 - e^{4\pi\mu}} \int_{0}^{2\pi}\int_{0}^{\min(t,\tau)} e^{2\mu(\tau - \sigma)} \dif\sigma (k^{2}q(\tau) + \mu^{2} - \alpha^{2})v(\tau) \dif\tau \\
           &  + \frac{e^{4\pi\mu}}{1 - e^{4\pi\mu}} \int_{0}^{t}\int_{\tau}^{t} e^{2\mu(\tau - \sigma)} \dif\sigma (k^{2}q(\tau) + \mu^{2} - \alpha^{2})v(\tau) \dif\tau \\
         = &\ 1 + \frac{1}{1 - e^{4\pi\mu}} \int_{0}^{2\pi} \frac{1}{2\mu}(e^{2\mu\tau} - e^{2\mu\max(\tau - t, 0)})(k^{2}q(\tau) + \mu^{2} - \alpha^{2})v(\tau) \dif\tau \\
           &  + \frac{e^{4\pi\mu}}{1 - e^{4\pi\mu}} \int_{0}^{t} \frac{1}{2\mu} (1 - e^{2\mu(\tau - t)})(k^{2}q(\tau) + \mu^{2} - \alpha^{2})v(\tau) \dif\tau \\
         = &\ 1 + \frac{1}{1 - e^{4\pi\mu}} \int_{0}^{2\pi} \frac{1}{2\mu}(e^{2\mu\tau} - e^{2\mu\max(\tau - t, 0)})(k^{2}q(\tau) + \mu^{2} - \alpha^{2})v(\tau) \dif\tau \\
           &  + \frac{e^{4\pi\mu}}{1 - e^{4\pi\mu}} \int_{0}^{2\pi} \frac{1}{2\mu}(e^{2\mu\max(\tau - t, 0)} - e^{2\mu(\tau - t)})(k^{2}q(\tau) + \mu^{2} - \alpha^{2})v(\tau) \dif\tau \\
         = &\ 1 + \int_{0}^{2\pi} \frac{1}{2\mu}(\frac{e^{2\mu\tau} - e^{4\pi\mu}e^{2\mu(\tau - t)}}{1 - e^{4\pi\mu}} - e^{2\mu\max(\tau - t, 0)})(k^{2}q(\tau) + \mu^{2} - \alpha^{2})v(\tau) \dif\tau.
\end{align*}
We define an integral operator $L$ by
\begin{equation*}
    (Lv)(t) \coloneqq \int_{0}^{2\pi} l(t, \tau) v(\tau) \dif\tau,
\end{equation*}
where
\begin{equation*}
    l(t, \tau) = \frac{1}{2}\left(\frac{e^{2\mu\tau} - e^{2\mu(2\pi - t + \tau)}}{1 - e^{4\pi\mu}} - e^{2\mu\max(\tau - t, 0)}\right)(k^{2}q(\tau) + \mu^{2} - \alpha^{2}).
\end{equation*}
Then we obtain the operator equation
\begin{equation*}
    v(t) - \frac{1}{\mu}(Lv)(t) = 1.
\end{equation*}
We observe that $\vert l(t, \tau)\vert \leqslant Ck^{2}$. 
 By the Neumann
series, we can show that
\begin{equation*}
    v(t) = \sum_{m = 0}^{\infty} \left(\frac{1}{\mu}L\right)^{m} = 1 + \frac{1}{\mu}\int_{0}^{2\pi} l(t, \tau) \dif\tau + O(k^{4}).
\end{equation*}

In the second case, if $v(d) = 0$, then the operator equation becomes $v(t) - \frac{1}{\mu}(Lv)(t) =
0$. Taking the norm of both sides of the equation $\mu v(t) = (Lv)(t)$, we find that $\mu \leqslant \Vert L\Vert
\leqslant Ck^{2}$. Sufficiently small $k$ guarantees that this case can not happen.

It remains to estimate the derivative $v'(t)$. We first compute $\partial_{t}
\frac{1}{\mu}(Lv)(t)$.
\begin{align*}
\partial_{t}\frac{1}{\mu} (Lv)(t) = &\ \partial_{t} \frac{1}{2\mu} \int_{0}^{2\pi} \frac{e^{2\mu\tau} - e^{2\mu(2\pi-t+\tau)}}{1 - e^{4\pi\mu}}(k^{2}q(\tau)+\mu^{2}-\alpha^{2})v(\tau) \dif\tau \\
                                    & - \partial_{t} \frac{1}{2\mu} \int_{0}^{t} (k^{2}q(\tau)+\mu^{2}-\alpha^{2})v(\tau) \dif\tau \\
                                    & - \partial_{t} \frac{1}{2\mu} \int_{t}^{2\pi} e^{2\mu(\tau-t)}(k^{2}q(\tau)+\mu^{2}-\alpha^{2})v(\tau)\dif\tau \\
                                  = &\ \int_{0}^{2\pi} \frac{e^{2\mu(2\pi-t+\tau)}}{1 - e^{4\pi\mu}}(k^{2}q(\tau)+\mu^{2}-\alpha^{2})v(\tau) \dif\tau \\
                                    & + \int_{t}^{2\pi} e^{2\mu(\tau-t)}(k^{2}q(\tau)+\mu^{2}-\alpha^{2})v(\tau) \dif\tau.
\end{align*}

It follows that
\begin{align*}
        \frac{\partial_{t}v(t)}{v(t)} = \frac{\partial_{t}\frac{1}{\mu}(Lv)(t)}{v(t)} = & \frac{1}{v(t)}\int_{0}^{2\pi} \frac{e^{2\mu(2\pi-t+\tau)}}{1 - e^{4\pi\mu}}(k^{2}q(\tau)+\mu^{2}-\alpha^{2})v(\tau) \dif\tau \\
                                        & +  \frac{1}{v(t)} \int_{t}^{2\pi} e^{2\mu(\tau-t)}(k^{2}q(\tau)+\mu^{2}-\alpha^{2})v(\tau) \dif\tau.
\end{align*}
Hence
\begin{equation*}
    \bigg\Vert \frac{\partial_{t}v(t)}{v(t)}\bigg\Vert_{L^{\infty}} \leqslant Ck^{2}.
\end{equation*}

We can use the argument in Step 3 to estimate $v_{n, 1}(d)$ and $v_{n, 1}'(d)$. We have that
\begin{equation*}
    \bigg\vert \frac{v_{n,1}'(d)}{v_{n,1}(d)}\bigg\vert \leqslant \bigg\Vert \frac{v_{n,1}'(t)}{v_{n,1}(t)}\bigg\Vert_{L^{\infty}} \leqslant Ck^{2}.
\end{equation*}
Hence $\big\vert \frac{v_{n,1}'(d)}{v_{n,1}(d)}\big\vert$ tends to $0$ as $k \to +0$ for fixed $n \neq 0$.
\end{Proof}

Now we study the asymptotic behavior for $n = 0$.

\begin{lemma}\label{lem:C3} It holds that $\theta_0=O(k)$ and 
 $\left\vert \frac{v_{0, 1}'(d)}{v_{0, 1}(d)}\right\vert=o(1)$  as $k \to +0$.
\end{lemma}

\begin{Proof}
(i)	It follows from Lemma \ref{lem:allEtan} that
	\begin{equation*}
		\lambda_{1} = \frac{\eta_{0} + \sqrt{\eta_{0}^{2} - 4}}{2},\ \lambda_{2} = \frac{\eta_{0} - \sqrt{\eta_{0}^{2} - 4}}{2}.
	\end{equation*}
	are mutually conjugate. From Theorem \ref{thm:discriminant} and the analyticity of
	$\eta_{0}(k)$, we know that $\eta_{0} \to 2$ as $k \to +0$ and $\theta_{0} \sim k$
	where $e^{i2\pi \theta_{0}} = \lambda_{1}$. 
	
(ii) For notational convenience we write $v=v_{0,1}$ to drop the dependence 
	on the subscripts. We shall present a sketch of the proof, following the argument 
    in Step 3 of Lemma \ref{lem:unibddk}. 

According to the downward modes in case (c) in Definition \ref{def:modes}, we 
assume that $u = e^{-i\theta_{0} t}v(t)$, where $v(0) = v(2\pi)$ and $v'(0) = v'(2\pi)$.
Substituting the product presentation of $u$ into $u''(t) + k^{2}(q(t) - \sin^{2}\theta)u(t) = 0$,
we have that
\begin{equation*}
    v''(t) - 2i\theta_{0} v'(t) + (k^{2}(q(t) - \sin^{2}\theta) - \theta_{0}^{2})v(t) = 0.
\end{equation*}
Then $(e^{-2i\theta_{0} t}v'(t))' = -e^{2i\theta_{0} t}(k^{2}(q(t) - \sin^{2}\theta)
- \theta_{0}^{2})v(t)$. Hence
\begin{equation*}
    e^{-2i\theta_{0}t}v'(t) = v'(0) - \int_{0}^{t} e^{-2i\theta_{0}\tau}(k^{2}(q(\tau) - \sin^{2}\theta) - \theta_{0}^{2})v(\tau) \dif\tau.
\end{equation*}
By straightforward computation, we get
\begin{align*}
    v'(t) =\ & \frac{1}{1 - e^{-4\pi i \theta_{0}}} \int_{t}^{2\pi} e^{-2i\theta_{0}(\tau - t)}(k^{2}(q(\tau) - \sin^{2}\theta) - \theta_{0}^{2})v(\tau) \dif\tau \\
           &+\frac{e^{-4\pi i \theta_{0}}}{1 - e^{-4\pi i \theta_{0}}} \int_{0}^{t} e^{-2i\theta_{0}(\tau - t)}(k^{2}(q(\tau) - \sin^{2}\theta) - \theta_{0}^{2})v(\tau) \dif\tau.
\end{align*}
Then we have
\begin{align*}
    v(t) = &\ 1 + \int_{0}^{t}v'(\sigma) \dif\sigma \\
         = &\ 1 + \frac{1}{2i\theta_{0}} \int_{0}^{2\pi} \left(e^{2i\theta_{0}\min(t - \tau, 0)} + \frac{e^{2i\theta_{0}(t - \tau)} - e^{2i\theta_{0}(2\pi - \tau)}}{e^{4\pi i \theta_{0}} - 1}\right)(k^{2}(q(\tau) - \sin^{2}\theta) - \theta_{0}^{2})v(\tau) \dif\tau.
\end{align*}
We define the integral operator $L$ by
\begin{equation*}
    (Lv)(t) \coloneqq \int_{0}^{2\pi} l(t, \tau) v(\tau) \dif\tau,
\end{equation*}
where
\begin{equation*}
    l(t, \tau) = \frac{1}{2i}\left(e^{2i\theta_{0}\min(t - \tau, 0)} + \frac{e^{2i\theta_{0}(t - \tau)} - e^{2i\theta_{0}(2\pi - \tau)}}{e^{4\pi i \theta_{0}} - 1}\right)(k^{2}(q(t) - \sin^{2}\theta) - \theta_{0}^{2}).
\end{equation*}
Then we obtain the operator equation
\begin{equation*}
    v(t) - \frac{1}{\theta_{0}}(Lv)(t) = 1.
\end{equation*}
We observe that $\vert l(t, \tau)\vert \leqslant Ck^{2}$. By a Neumann series
argument, we have that
\begin{equation*}
    v(t) = \sum_{m = 0}^{\infty} \left(\frac{1}{\theta_{0}}L\right)^{m} = 1 + \frac{1}{\theta_{0}}\int_{0}^{2\pi} l(t, \tau) \dif\tau + O(k^{2}).
\end{equation*}
By arguments similar to Step 3 of the proof of Lemma \ref{lem:unibddk}, We can 
estimate $v'(t)$ and obtain that $\left\vert \frac{v_{0, 1}'(d)}{v_{0, 1}(d)}\right\vert$
tends to $0$ as $k \to +0$.
\end{Proof}

\end{appendices}


\end{document}

%% file: figs/bvp.tex
	\begin{tikzpicture}
		
		\draw[line width=1.75cm,color=red!20] (0,0.875) -- (5,0.875);
		\foreach \x in {1,2,3}{
			\draw[line width=0.25cm,color=red!50] (0,-0.125+\x*0.5) -- (5,-0.125+\x*0.5);
		}
		\draw[line width=1.75cm,color=red!20] (0,3.625) -- (5,3.625);
		
		\draw[line width=1cm,color=cyan!10] (0,2.25) -- (5,2.25);
		
		\draw (0,1.75) -- (5,1.75);
		\draw (0,2.75) -- (5,2.75);
		
		\draw[-{stealth},blue] (0.5,3.8) -- (1.5, 2.6);
		\draw[blue] (1.5,2.8) node [above] {\footnotesize$u^{\mathrm{in}}$};
		
		\draw[-{stealth},blue] (2.5,2.3) -- (3.5,3.5);
		\draw[blue] (3.5,3.5) node [above] {\footnotesize$u^{\mathrm{sc}}$};
		
		\draw[-{stealth},blue] (2,2) -- (3,1);
		\draw[blue] (3,1) node [below] {\footnotesize$u^{\mathrm{tr}}$};
		
		\draw (5,2.75) node [right] {\tiny$\Gamma_{d}$};
		\draw (4,3.625) node [right] {$\Omega_{d}^{+}$};
		\draw (4,2.5) node [right] {$\Omega$};
		\draw (5,1.75) node [right] {\tiny$\Gamma_{b}$};
		\draw (4,0.875) node [right] {$\Omega_{b}^{-}$};
	\end{tikzpicture}

%% file: figs/bvpc.tex
	\begin{tikzpicture}
		\draw[line width=1.75cm,color=red!20] (0,0.875) -- (5,0.875);
		\foreach \x in {1,2,3}{
			\draw[line width=0.25cm,color=red!50] (0,-0.125+\x*0.5) -- (5,-0.125+\x*0.5);
		}
		\draw[line width=1.75cm,color=red!20] (0,3.625) -- (5,3.625);
        \draw[dashed] (0,0) -- (0,4.5);
        \draw[dashed] (5,0) -- (5,4.5);
		
		\draw[line width=1cm,color=cyan!10] (0,2.25) -- (5,2.25);
		\draw (0,1.75) -- (5,1.75);
		\draw (0,2.75) -- (5,2.75);
		
		\draw (5,2.75) node [right] {\tiny$\tilde{\Gamma}_{d}$};
		\draw (4,2.3) node [right] {$C$};
		\draw (5,1.75) node [right] {\tiny$\tilde{\Gamma}_{b}$};
	\end{tikzpicture}

%% file: figs/discriminant.tex
	\begin{tikzpicture}
		\draw[thick,->] (-5,0) -- (5,0);
		\draw[thick,->] (-4,-2) -- (-4,2);
		\draw[red] (-5,1) -- (5,1);
		\draw[red] (-5,-1) -- (5,-1);

		\draw[blue, thick] plot[smooth] coordinates {(-4,1.5) (-3.5,1) (-2.5,0) (-1.5,-1) (-0.5,-1) (1.5,1) (2.5,1) (3.5,0.5)};
		\draw[dashed] (-3.5,0) -- (-3.5,1);
		\draw[dashed] (-1.5,0) -- (-1.5,-1);
		\draw[dashed] (-0.5,0) -- (-0.5,-1);
		\draw[dashed] (1.5,0) -- (1.5,1);
		\draw[dashed] (2.5,0) -- (2.5,1);

		\draw (-4,2) node [above] {$\eta(\lambda)$};
		\draw (-3.5,0) node [below] {$\lambda_{0}$};
		\draw (-1.5,0) node [above] {$\mu_{0}$};
		\draw (-0.5,0) node [above] {$\mu_{1}$};
		\draw (1.5,0) node [below] {$\lambda_{1}$};
		\draw (2.5,0) node [below] {$\lambda_{2}$};
		\draw[red] (5,1) node [right] {$2$};
		\draw[red] (5,-1) node [right] {$-2$};
	\end{tikzpicture}